\documentclass[12pt]{article}
\usepackage{amsmath,amssymb,amsthm,amsfonts}
\usepackage[numeric]{amsrefs}
\usepackage{hyperref}

\newtheorem{thm}{Theorem}[section]
\newtheorem{cor}[thm]{Corollary}
\newtheorem{lem}[thm]{Lemma}
\newtheorem{prop}[thm]{Proposition}

\newcommand{\thmref}[1]{Theorem~\ref{#1}}
\newcommand{\lemref}[1]{lemma~\ref{#1}}
\newcommand{\propref}[1]{proposition~\ref{#1}}

\newcommand{\secref}[1]{section~\ref{#1}}

\newcommand{\bx}{\hfill$\square$\vspace{.6cm}}

\DeclareMathOperator{\sg}{sgn}

\numberwithin{equation}{section}

\newcommand{\gobble}[1]{}
  \newcommand{\rangeref}[2]{%
    \ref{#1}--\afterassignment\gobble\fam 0\ref{#2}%
  }

\renewcommand\a{\alpha}         
\renewcommand\b{\beta}
\newcommand\g{\gamma}
\renewcommand\d{\delta}
\newcommand\e{\varepsilon}
\renewcommand\l{\lambda}

\newcommand\G{\Gamma}
\newcommand\f{\frac}
\newcommand{\N}{{\mathbb{N}}}
\newcommand{\Z}{{\mathbb{Z}}}
\newcommand{\R}{{\mathbb{R}}}

\newcommand{\C}{{\mathbb{C}}}

\newcommand{\Q}{{\mathbb{Q}}}

\newcommand\im{\mbox{Im~}}
\renewcommand\Re{\text{Re~}}

\newcommand{\ttwo}[4]{\left(\begin{array}{cc}
{#1} & {#2} \\ {#3} & {#4} \end{array} \right)}

\renewcommand\i{^{-1}}

\renewcommand\({\left(}         
\renewcommand\){\right)}

\begin{document}

\author{Stephen D. Miller\thanks{Partially supported
by National Science Foundation grants DMS-0122799, DMS-0301172,
and an Alfred P. Sloan Foundation Fellowship. }}

\title{Cancellation in additively twisted sums on~$GL(n)$}
\date{August 29, 2004}

\maketitle

\begin{abstract}

In a previous paper with Schmid \cite{regularity} we considered
the regularity of automorphic distributions for $GL(2,\R)$, and
its connections to other topics in number theory and analysis. In
this paper we turn to the higher rank setting, establishing the
nontrivial bound $\sum_{n\le T}a_{n}\, e^{\,2\,\pi\, i\, n \,\a} =
O_\e(T^{\,3/4+\e})$, uniformly in $\a\in\R$, for $a_n$ the
coefficients of the $L$-function  of a cusp form on
$GL(3,\Z)\backslash GL(3,\R)$.  We also derive an equivalence
(\thmref{perthm}) between analogous cancellation statements for
cusp forms on $GL(n,\R)$, and the sizes of certain period
integrals. These in turn imply estimates for the second moment of
cusp form $L$-functions.

\end{abstract}

\section{Introduction}\label{intro}

 Consider a sequence $a_n$ of
arithmetic quantities of order 1, and the sums of their twists by
additive characters
\begin{equation}\label{sdef}
S(T,\a) \ \ = \ \ \sum_{n\,=\,1}^T \ a_n \  e^{\,2\, \pi \, i \, n
\, \alpha} \ \ , \ \ \ \ \a \ \in \ \R\,. \end{equation} In this
paper we shall be concerned with obtaining estimates for
$S(T,\alpha)$ which are uniform in $\a$.  This problem was
considered already by Hardy and Littlewood in 1914 \cite{HL}, and
is well understood when the $a_n$ are the normalized Fourier
coefficients of a modular or Maass form on the upper half plane,
i.e.\! automorphic forms on $GL(2,\R)$ (see \cite{regularity}). In
the case of cusp forms (which is simpler to state), one has the
estimate
\begin{equation}\label{best}
    S(T,\alpha) \ \ = \ \ O_\e(T^{1/2+\e}) \ , \  \ \ \text{for
    any~}\, \e \, > \, 0\,,
\end{equation}
uniformly in $\a$.\footnote{The equivalent notations $A=O(B)$ and
$A\ll B$ signify that there exists a positive constant $C$ such
that $|A| \, \le \, C \, B$.  We write $A=O_\e(B)$ or $A\ll_\e B$
when this holds for all positive values of a parameter $\e$, but
with a constant $C=C_\e$ that potentially depends on $\e$.} This
can be seen in a variety of ways, perhaps most naturally in terms
of the boundedness of cusp forms (e.g. (\ref{addedholom})). The
exponent of $T^{1/2}$ is best-possible, as can be seen by
estimating the $L^2$-norm of the trigonometric polynomial
$S(T,\a)$
\begin{equation}\label{l2norm}
    \int_0^1 \, |S(T,\a)|^2\,d\a \ \ = \ \
    \sum_{n \, \le \, T } |a_n|^2\,,
\end{equation}
which should be of order $T$ if the $a_n$ are of order 1.
 See, for
example,
\cite{erdos,nara2,hafner,murty,chowla,HL,walfisz,regularity} for
background on techniques used to bound $S(T,\alpha)$.  A folklore
conjecture asserts that the estimate (\ref{best}) holds for the
Fourier coefficients of any cusp form, on any  group.  The purpose
of this paper is to provide a non-trivial, uniform estimate for
$S(T,\a)$ beyond the classical case of $GL(2,\R)$
(\thmref{mainthm} below).

 Such sums have long
been connected to important questions in analytic theory.  Most
notably, Titchmarsh's method \cite[p.\! 165]{Titch} derives from
(\ref{best}) the correct order of magnitude
\begin{equation}\label{magn}
    \int_{-T}^T\,\left|
   L(1/2+it) \right|^2 \, dt \ \ = \ \ O_\e(T^{1+\e}) \ \ , \ \ \
   \e \, > \, 0
\end{equation}
for the second moment of the $L$-function $L(s) =
\sum_{n=1}^\infty a_n n^{-s}$ formed from the cusp form's
coefficients $a_n$.  We include a proof of this in \thmref{implic}
below. A slight variant for the coefficients of arbitrary
automorphic forms on $GL(m)$, $m\ge 2$ implies estimates for the
higher moments of $L(s)$ as well:
\begin{equation}\label{magnk}
    \int_{-T}^T\,\left|
   L(1/2+it) \right|^{\,2k} \, dt \ \ = \ \ O_\e(T^{1+\e}) \ \ , \ \ \
   \e \, > \, 0 \ \ , \, \ \ \text{for all~}k\, \ge \, 1\,.
\end{equation}  The latter is equivalent to the
generalized Lindel\"of conjecture in the $t$-aspect, which states
that
\begin{equation}\label{lindlt}
    L(1/2+it) \ \ = \ \ O_\e(t^{\e}) \ \ , \ \ \ \e \, > \, 0
    \, .
\end{equation}
It is widely believed that it is just as difficult to obtain, for
example, the correct order of magnitude (\ref{magnk}) of the
second moment of the standard $L$-function of an automorphic form
on $GL(3)$, as it is to obtain the correct order of magnitude for
the sixth power moment of the Riemann $\zeta$-function. This has
long been a major challenge in analytic number theory.

 Thus the estimate (\ref{best}), not for modular forms but
for automorphic cusp forms of higher rank, is evidently a very
difficult one to obtain.  It is not surprising that good bounds
for $S(T,\a)$, $\a \in \Q$, can be obtained;
 a classical result
of Landau (see \cite[Chapter 12]{Titch}) gives bounds of the form
$O(T^\Theta)$, where $\Theta<1$.  This cancellation is closely
related to the analytic continuation and functional equation of
the multiplicatively twisted $L$-functions $L_\chi(s) \ \ = \ \
\sum_{n=1}^\infty a_n\, \chi(n)\,n^{-s}$, where $\chi$ is a
Dirichlet character. The main challenge is to also provide
estimates when $\a$ is irrational, and uniform ones at that.

In this paper we will deduce such bounds  by approximating $\a$ by
rational numbers.  As an illustration, consider a holomorphic cusp
form $f(z)=\sum_{n\ge 1}a_n n^{(k-1)/2}e(nz)$ of weight $k$ for
$SL(2,\Z)$.  (The $a_n$ are the coefficients of the standard
$L$-function of $f$.) Such a form satisfies the uniform bound
$|f(z)|=O((\im z)^{-k/2})$ for all $z=\a+I/T$ in the upper half
plane, and so in particular for any $\a\in \R$
\begin{equation}\label{addedholom}
    \sum_{n\,=\,1}^\infty  \left[\,a_n \,e(n\a) \, \right]\,n^{(k-1)/2}\,
    e^{-2\,\pi\,n/T}
    \ \ = \ \ O(T^{k/2}) \ , \ \ \ \text{as~}t\, \rightarrow \, 0\,.
\end{equation}
This is essentially a smoothed form of (\ref{best}).  The bound
for $f(z)$, in turn, comes from the modularity of $f$ under a
suitably chosen matrix $\g\in SL(2,\Z)$ which maps $z$ to a point
$\g z$ in a fixed fundamental domain, on which $f(z)$ is bounded.
Finding such a $\g$ is a diophantine problem.  In particular, if
$\a\in \Q$, $\g$ can be chosen so that $\g z$ is very close to the
cusp, where $f(z)$ in fact decays rapidly; this partly explains
the remark of the previous paragraph.
 In order to generalize this argument to non-holomorphic cusp forms such as Maass forms, or to
 automorphic forms on $GL(3)$, we will  use a Voronoi-style
summation formula (\secref{vsec}) to give bounds on smoothed sums
analogous to (\ref{addedholom}).

Our main result (\thmref{mainthm} below) is the non-trivial
uniform estimate for $S(T,\a)$ of $O_\e(T^{3/4+\e})$ when the
$a_n$ are the Fourier coefficients of an automorphic form $\Phi$
on $GL(3,\Z)\backslash GL(3,\R)$.\footnote{The techniques of this
paper and \cite{voronoi} apply to the general congruence subgroup
of $GL(3,\Z)$, but the coefficients $a_{q,n}$ are no longer
uniquely determined by the $L$-function data.} These, which we
shall denote $a_{q,n}$, are naturally indexed by two integral
parameters, the first of which we will hold fixed and not attempt
to measure the dependence of (though it is of course possible to
do so). The coefficients $a_{n,q}$ are actually the Fourier
coefficients of the form $\tilde\Phi$ contragredient to $\Phi$, so
there is no loss of generality in fixing the first index instead
of the second.
 Thanks to the Rankin-Selberg
theory (\cite{Jac-Sha}; see also \cite[\S 2]{bumprs} and \cite[\S
5]{jacquetindia}),
 we know the $a_{q,n}$ obey the Ramanujan
conjecture on average.  More precisely, the Rankin-Selberg
$L$-function $$L(s,\Phi\otimes\tilde\Phi)  \ \ = \ \ \sum_{n,\,q\,
\ge \, 1} \ |a_{q,n}|^2\ (nq^2)^{-s}\, ,$$ initially convergent
for $\Re{s}$ large, has a meromorphic continuation to $\C$ with
only a simple pole at $s=1$; this translates into the estimates
\begin{equation}\label{rspole}\gathered
    \sum_{nq^2\, \le \, T}  \, |a_{q,n}|^2 \ \ = \ \ O(T)\, , \\
    \sum_{n \, \le \, T}  \, |a_{q,n}|^2 \ \ = \ \ O(q^2\, T)\, ,
    \endgathered
\end{equation}
and
\begin{equation}\label{ramave}
    \sum_{n\,  \le \,  T} \,  |a_{q,n}|  \ \ =  \ \ O(   q \,  T  )
\end{equation}
by Cauchy-Schwartz.  Thus the trivial estimate for $S(T,\alpha)$
-- obtained by taking the absolute value of each term in
(\ref{sdef}) -- is $O(T)$ for any fixed $q$. Our main result is
the following improvement, which goes halfway between the trivial
bound and the best-possible bound of $O(T^{1/2})$.

\begin{thm}\label{mainthm}
Let $a_{q,n}$ denote the Fourier coefficients of a cusp form on
$GL(3,\Z)\backslash GL(3,\R)$.  Then for any $\e>0$
\begin{equation}\label{mainbound}
\sum_{n \,=\,1}^T \,a_{q,n}  \ e^{\,2 \, \pi \,   i \,  n  \, \a}
\ \ = \ \ O_\e(T^{\, 3/4+\e}) \ ,\ \ \ \ \hbox{uniformly in~}\a
\in \R\, ,
\end{equation}
with the implied constant depending on  $q$, $\e$, and the cusp
form.
\end{thm}
\noindent Through partial summation, this statement implies a
bound for the analogous smoothed sums (actually \thmref{mainthm}
is derived from a similar statement -- see (\ref{punch1})):
\begin{cor}\label{smoothedthm}
Let $a_{q,n}$ denote the Fourier coefficients of a cusp form on
$GL(3,\Z)\backslash GL(3,\R)$,  and $\phi$ be any Schwartz
function.  Then
 for any $\e>0$ \begin{equation}\label{smoothedbd}
    \sum_{n\neq 0} \, a_{q,n} \, e^{\,2\,\pi\,i\,n\,\a}\,\phi\(\f nT\)
    \ \ = \ \ O_\e(T^{\,3/4+\e})\, ,
\end{equation}
where the implied constant depends on $\e$, $q$, $\phi$, and the
cusp form.
\end{cor}

As we mentioned before, the folklore conjecture that
$S(T,\alpha)=O_\e(T^{1/2+\e})$ implies the correct order of
magnitude for the second moment (\ref{magn}).  Weaker estimates on
$S(T,\a)$ still give cancellation bounds via the classical method
of Titchmarsh alluded to above.  Though the following theorem
appears to be well-known to experts, we have been unable to locate
a suitable statement in the literature, and so have chosen to
include a proof in \secref{moments}.

\begin{thm}\label{implic}
Let $L(s)=\sum_{n=1}^\infty \,a_n\,n^{-s}$ be the $L$-function of
a cusp form on $GL(m)$ over $\Q$, other than the Riemann
$\zeta$-function.  Suppose that
\begin{equation}\label{implic1}
    S(T,\a) \ \ = \ \ \sum_{n\,=\,1}^T\,a_n\,e(n\,\a) \ \ = \ \
    O_\e(T^{\b+\e}) \ \ , \ \ \ \ \text{uniformly in} \ \a \,,
\end{equation}
for some $\b\ge \f12$ and any $\e>0$.  Then the second moment of
$L(s)$ satisfies the bound
\begin{equation}\label{implic2}
    \int_{-T}^T\,|L({\scriptstyle \f 12}+it)|^2\,dt \ \ = \ \
    O_\e(T^{\,1 \,+\: (2\b-1)\,m\,+\,\e})\ \ , \ \ \ \ \text{for
    any} \ \e \, > \, 0\,.
\end{equation}
\end{thm}
Some brief remarks are in order.  First, the omission of
$\zeta(s)$ is made for a technical reason; besides the fact that
the precise asymptotics of the second moment of $\zeta(\f 12 +
it)$ have long been known (see \cite[\S7]{Titch}), \thmref{implic}
requires some modification for $L$-functions which have poles.
Such an adjustment can be used to study the $2k$-th moment
(\ref{magnk}), though we shall not pursue this here. Second, we
have chosen to state \thmref{implic} for cusp form $L$-functions
on $GL(m)$ over $\Q$ because they and their products are believed
to account for the totality of $L$-functions. Titchmarsh's method
could equally be used to derive results for other classes, such as
the Selberg class \cite{selclass}.  Thirdly, \thmref{implic} is of
interest only for $\b \le \f 34 - \f{1}{2m}$ (and hence $m > 2$),
because the second moment (\ref{implic2}) can always be bounded by
$O_\e(T^{m/2+\e})$ using the approximate functional equation (see
\cite[p. 31]{nara}).

Though \thmref{mainthm} is the first nontrivial bound for
$S(T,\a)$ on  $GL(m)$, $m
> 2$, it still falls far short of improving any  estimates on
the critical values of a $GL(3)$ $L$-function. Our obstacle to
sharpening the estimate of \thmref{mainthm} is the appearance of
Kloosterman sums in formula (\ref{v}), which we bound  in
\secref{pfsec} only by their absolute value (Weil's bound). Future
improvements would necessarily obtain cancellation in sums of
products of the $a_n$ with Kloosterman sums.  We are unable to
prove any interesting statements for $GL(m)$, $m > 3$, but there
obtaining cancellation in sums of $a_n$ times {\em
hyper-}Kloosterman sums could in principle be used to attack the
second moment.  Though this appears no easier, it is perhaps of
interest that the moment problem is connected to exponential sums
in this fashion.

Sections 2-5 of this paper contain the proof of \thmref{mainthm};
in section 6 we turn to the proof of \thmref{implic}.  Finally, in
section 7  we give an equivalence between bounds on $S(T,\a)$ for
cusp forms, and the sizes of certain period integrals studied by
Jacquet, Piatetski-Shapiro, and Shalika in their construction of
the standard $L$-function on $GL(m)$. In particular, the
equivalence given by \thmref{perthm}b for the optimal case of
$S(T,\a)=O_\e(T^{1/2+\e})$ can be viewed as a condition on an {\em
individual} cusp form which, together with \thmref{implic},
implies the correct order of magnitude for the second moment of
its standard $L$-function. (Again, a modification for non-cusp
form $L$-functions can be used to discuss higher moments and the
full Lindel\"of conjecture in the $t$-aspect.)

Our interest in this problem originated in joint work with
Wilfried Schmid on questions regarding the H\"older regularity of
the boundary distributions associated to cusp forms on $GL(3,\R)$
(see \cite{regularity} for a survey on the case of $GL(2,\R)$).
Theorem~\ref{mainthm} can be used to give the following
non-trivial estimate:

\begin{cor}\label{cor1}  Let
$$\tau_{x,q}(x) \ \ = \ \ \sum_{n\, \neq \, 0}\, c_{n,q}\, e(n\,x)$$
denote the abelian Fourier components of the of the boundary value
distribution of an automorphic cusp     form on
$GL(3,\Z)\backslash GL(3,\R)$, where $$c_{n,q}  \ \ = \  \ a_{n,q}
\, |n|^{-\l_1} \, |q|^{\l_3} \, \sg(n)^{\d_1} \, \sg(q)^{\d_3}$$
(the $\l_j$ and $\d_j$ are representation-theoretic parameters as
in the next section -- see \cite[\S7]{voronoi} for details). Then
$\tau_{x,q}$ lies in the H\"older class $C^{\,<\Re{\l_1}-3/4}$.
\end{cor}

\noindent The definitions of the H\"older classes as well as the
proof of this corollary can also be found in \cite[\S
3]{regularity}. The above folklore cancellation conjecture that
$S(T,\a)=O_\e(T^{1/2+\e})$ can be restated in terms of boundary
regularity as the assertion that $\tau_{x,q}\, \in \,
C^{\,<\Re{\l_1}-1/2}$.  Interestingly, the techniques from partial
differential equations and representation theory used in
\cite{flato} -- which obtain an essentially sharp estimate for
$GL(2,\R)$ -- seem to only recover  a very weak bound for
$GL(m,\R)$, $m\ge 3$. This is consistent with the expected overall
difficulty of (\ref{magn}), which is a consequence of
$S(T,\a)=O_\e(T^{1/2+\e})$.

Finally, a remark is in order about the coefficients of
noncuspidal automorphic forms.   For example, the early papers of
\cite{erdos,chowla,HL,walfisz}  studied the Fourier coefficients
of Eisenstein series on the upper half plane, notably the divisor
function $d(n)$.
 There is an Eisenstein series on
$GL(3,\Z)\backslash GL(3,\R)$, whose standard $L$-function is
$\zeta(s)^3$, with Fourier coefficients $a_{1,n}$  equal to the
triple divisor function
\begin{equation}\label{tripdiv}
d_3(n) \ \ := \ \ \# \, \{\,a \,b \,c \, = \,  n  \  |  \ a,b,c \,
\in \, \N\,\}\,.
\end{equation}
 The method used in this paper can
be extended to study additively twisted sums of $d_3(n)$ as well,
although there is of course no nontrivial bound for even $S(T,0)$
here because $d_3(n)>0$.  One must settle for almost-everywhere
bounds, which could not possibly be uniform (or else by continuity
they would extend everywhere).  Strong nonuniform results,
however, can be obtained via Carleson's theorem on Fourier series
(see \cite{murty}) by using only the fact that
$d_3(n)=O_\e(n^\e)$, and nothing about automorphy.  In any event,
for many applications  -- such as in studying moments -- a more
useful form of (\ref{best}) would uniformly bound the difference
between $S(T,\a)$ and a main term.

\subsection*{Acknowledgements}  It is a pleasure to thank J. Cogdell, W. Duke,
J. Friedlander, H. Iwaniec, E. Lapid, M. R. Murty, R. Narasimhan,
P. Sarnak, W. Schmid, and K. Soundararajan for their helpful
comments and discussions.

\section{Voronoi Summation for $GL(3)$}\label{vsec}

Our main tool is the $GL(3)$ analog of the Voronoi summation
formula, recently proven in \cite{voronoi} (see also \cite{expos}
for a less technical exposition).  We will now give a brief
summary of the formula, referring the reader to \cite{voronoi} for
definitions not fully explained here, and for its connections to
the functional equations of twisted $L$-functions. Denote the
embedding parameters of a cusp form $\Phi$ on $GL(3,\Z)\backslash
GL(3,\R)$ by $(\l_1,\l_2,\l_3)\in\C^3$ and $(\d_1,\d_2,\d_3)\in
(\Z/2\Z)^3$. These are representation-theoretic parameters
connected to $\Phi$, and will be related to the functional
equation for the $L$-function $L(s,\Phi)$ below in (\ref{funeq}).
We may and shall make the normalizing assumption that
\begin{equation}\label{normassumpt}
\l_1 \, + \, \l_2 \, + \, \l_3 \ = \ 0  \ \ \  \ \hbox{~and~}  \ \
\ \ \d_1 \, + \,  \d_2 \, + \, \d_3 \ \equiv  \ 0 \!\! \pmod 2.
\end{equation}
 The summation formula which we are about to state involves
doubly-indexed Fourier coefficients  $a_{n,m}$  of a cusp form
$\Phi$ on $GL(3,\Z)\backslash GL(3,\R)$.  These are perhaps
simplest described in terms of the standard and contragredient
$L$-functions of $\Phi$,
\begin{equation}\label{L-functions}
    L(s,\Phi) \ \ = \ \ \sum_{n=1}^\infty \,a_{1,n}\,n^{-s} \ \ \
    \text{and} \ \ \ L(s,\tilde{\Phi})\ \
    = \ \ \sum_{n=1}^\infty \,a_{n,1}\,n^{-s} \, ,
\end{equation}
respectively.  If $\Phi$ is  a Hecke eigenform -- an assumption we
do not make, yet one which entails no loss of generality -- the
coefficients $a_{n,m}$ are
 eigenvalues of the Hecke operators $T_{n,m}$, and accordingly
satisfy certain recursion identities (for a full description see
\cite[\S 9]{Bump}).  In particular, when $\Phi$ is a Hecke
eigenform, the $a_{n,m}$ can be derived from the $a_{n,1}$ and
$a_{1,m}$ via the identity
\begin{equation}\label{anmrecurs}
    a_{n,m} \ \ = \ \ \sum_{d\,| \, (n,m)}\,\mu(d)\,a_{n/d,1} \,
    a_{1,m/d}\,,
\end{equation}
where $\mu(d)$ denotes the M\"obius $\mu$-function.

The following is the Voronoi-style summation formula for
automorphic forms on $GL(3,\Z)\backslash GL(3,\R)$.

\begin{thm}\label{vsf} (\cite{voronoi}) Suppose that $a_{n,m}$ are the Fourier
coefficients of a cuspidal $GL(3,\Z)$-automorphic representation
of $GL(3,\R)$, with embedding parameters $(\l_1,\l_2,\l_3)$ and
$(\d_1,\d_2,\d_3)$ as in (\ref{normassumpt}). Let $f$ be a
Schwartz function which vanishes to infinite order at the origin,
or more generally, a function on $\,\R - \{0\}\,$ such that
$\,(\sg x)^{\d_3}|x|^{-\l_3}f(x)$ is a Schwartz function. Then for
$T>0$, $(a,c)=1$, $c \neq 0$, $a \bar  a \equiv 1 \
(\operatorname{mod} c)$, and $q
> 0$,
\begin{equation}\label{v}
\sum_{n\neq 0} a_{q,n}\,e(-na/c)\,f(\f nT)\ =\
\sum_{d|cq}\left|\f{c}{d}\right| \sum_{n\neq 0} \f{a_{n,d}}{|n|}\,
S(q\bar{a},n;qc/d) \, F\!\(\f{nd^2}{c^3q}T\),
\end{equation}
where $\,S(n,m;c)\,=\,\sum_{x\in(\Z/c\Z)^*}
e\(\f{nx+m\bar{x}}{c}\)\,$ denotes the Kloosterman sum and, in
symbolic notation,
\begin{equation}\label{vvv}
F(t)=\int_{\R^3}f\!\(\f{x_1x_2x_3}{t}\)\ {\prod}_{j=1}^3 \((\sg
x_j)^{\d_j}\,|x_j|^{-\l_j}\,e(-x_j)\) dx_3 \, dx_2 \, dx_1 \, .
\end{equation}
This integral expression for $F$ converges when performed as
repeated integral in the indicated order -- i.e., with $x_3$
first, then $x_2$, then $x_1$ -- and provided
$\,\operatorname{Re}\l_1 > \operatorname{Re}\l_2 >
\operatorname{Re}\l_3$\,; it has meaning for arbitrary values of
$\,\l_1,\,\l_2,\,\l_3\in\C\,$ by analytic continuation.  (An
alternative description of the relation (\ref{vvv}) is given in
(\ref{MFandMf}) below.)
\end{thm}

The embedding parameters $(\l_1,\l_2,\l_3)$ obviously play an
important role in \thmref{vsf}, so it is worthwhile to describe
them in more detail. The parameter $\l_3$ may always be chosen to
have  the maximal real part among $\{\l_1,\l_2,\l_3\}$. We can
assume
\begin{equation}\label{orderoflambda}
    \Re{\l_1} \, , \,  \Re{\l_2}  \ \le \ \Re{\l_3}  \ , \ \ \
    \Re{\l_1} \, , \,  \Re{\l_2} \  < \ \hbox{$\f 12$}  \  , \ \,  \text{~and}
   \ \ \ \  \Re{\l_3} \ \ge \ 0\, .
\end{equation}
The second inequality requires some explanation. There are
essentially two types of representations of $GL(3,\R)$
corresponding to  cusp forms.  The first possibility is that
$\Phi$ comes from a fully induced principal series representation,
in which case $|\Re{\l_j}|<1/2$ and the $\l_j$'s may be freely
permuted; otherwise $\Phi$ is connected to a induced
representation of $GL(3,\R)$ constructed from the discrete series
$D_k$ of $GL(2,\R)$ (corresponding to weight $k \ge 2 $ modular
forms). In this latter situation we may and do chose the $\l_j$'s
to be written as
\begin{equation}\label{lambdared}
    \l_1 \ = \ -2\,i \,t   \ , \ \  \l_2 \ = \ -\,\f{k-1}{2} + i t \ , \ \
    \l_3 \  = \ \f{k-1}{2} + i t \ , \ \ t \in \R\,,
\end{equation}
with $\d_1\equiv \d_2+\d_3\equiv k \pmod 2$.
 In either case (\ref{orderoflambda}) certainly holds.
 This
bound of $1/2$ comes from knowledge of the unitary dual of
$GL(3,\R)$ and has a generalization to $GL(n,\R)$.  For
automorphic representations one can in fact do better, though this
is not necessary for our purposes (see \cite[Appendix 2]{Kim} for
the most recent improvements).

A more useful characterization of the relation between the
functions $f$ and $F$ in (\ref{vvv}) is provided by the (signed)
Mellin transforms.  For this we must split the functions $f$ and
$F$ into odd and even parts; the relation (\ref{vvv}) clearly
preserves parity.  If a function $g$ has parity $\eta \in \Z/2\Z$
(i.e. $g(-x)=(-1)^\eta g(x)$), then we define the signed Mellin
transform of $g$ to be
\begin{equation}\label{signedmelldef}
    M_{\eta \,} g(s)  \ \ = \ \ \int_\R \, g(x) \, |x|^{s-1} \,
    \sg(x)^\eta \, dx
\end{equation}
for values of $s$ where the integral is absolutely convergent, and
elsewhere by analytic continuation. When both $f$ and $F$ have
parity $\eta$, the relation (\ref{vvv}) can be succinctly
described by the
 formula
\begin{equation}\label{MFandMf}
    M_{\eta\,} F(s-1) \ \ = \ \ G_{\d_1+\eta}(s-\l_1) \,
    G_{\d_2+\eta}(s-\l_2) \, G_{\d_3+\eta}(s-\l_3) \ M_{\eta\,} f (1-s)\,
\end{equation}
(\cite[Theorem 1.18]{voronoi}).  Here the function
\begin{equation}\label{Gdcossin}
\aligned
    G_\d(s) \ \ &= \ \  \  & \ \ \  (2\pi)^{\, -s} \, \G(s) \,
    \left[ \, e\(\f s4 \) \,  +  \, (-1)^\d  \, e\(- \f s4 \) \, \right]
    \qquad   \\ \\
&= \ \ &\left\{
\begin{array}{ll}
    \, \ 2 \, (2\pi)^{-s} \, \G(s)\,\cos(\pi s/2)\, , & \d\equiv 0\pmod 2\, \\
   \\
   2\,i \, (2\pi)^{-s} \, \G(s)\,\sin(\pi s/2)\, , & \d\equiv 1\pmod 2\, \\
\end{array}%
\right.
\endaligned
\end{equation}
has only simple poles and simple zeroes, at $s \in (2\Z+\d)\cap \Z_{\le 0}$
and $s\in (2\Z+\d+1)\cap \Z_{>0}$, respectively. The functional
equation relating $L(s,\Phi)$ and $L(s,\tilde\Phi)$ can be cleanly
stated in terms of the $G_\d$ as
\begin{equation}\label{funeq}
    L(1-s,\tilde\Phi) \ \ = \ \ G_{\d_1}(s+\l_1)
                \, G_{\d_2}(s+\l_2) \,
                G_{\d_3}(s+\l_3) \, L(s,\Phi)\,;
\end{equation}
from this it is also possible to relate the $\l_j$ and $\d_j$ to
the
 $\G$-factors appearing in the usual form of the functional equation  (see
\cite[\S6]{voronoi}).

The functions $G_\d$ also arise in relating the Mellin and Fourier
transforms.  Suppose that $g$ is a Schwartz function of parity
$\eta$; then
\begin{equation}\label{ghatdef}
    \widehat{g}(r) \ \ = \ \ \int_\R \, g(x)\,e(-x r)\,dx
\end{equation}
 is also, and
\begin{equation}\label{mellhatg}
    M_{\eta\,}\widehat{g}(s) \ \ = \ \ (-1)^\eta \  G_\eta(s) \
    M_{\eta\,} g(1-s)
\end{equation}
(\cite[(4.58)]{inforder}).  The Fourier inversion formula is then
equivalent to the identity
\begin{equation}\label{grecip}
    G_\eta(s) \, G_\eta(1-s) \ \ = \ \ (-1)^\eta
\end{equation}
(\cite[(4.11)]{inforder}).

We end this section with a remark about the product of $G_\d$'s
occurring in (\ref{MFandMf}), namely that
\begin{equation}\label{holominsg1}
    G_{\d_1+\eta}(s-\l_1) \,
    G_{\d_2+\eta}(s-\l_2) \, G_{\d_3+\eta}(s-\l_3) \  \text{
    is holomorphic for~} \  \Re{s} \, \ge  \, \hbox{ $\frac 12$}\,.
\end{equation}
In light of our assumptions and discussion around
(\ref{orderoflambda}), the only possible poles must come from the
third factor, and even then only when
 $\Re{\l_3}\ge \f 12$, in which case we also assume (\ref{lambdared}).
 In this case the product of the last two
functions in (\ref{holominsg1}) in fact equals
$$ i^k \, (2\pi)^{\,1\,-\,2\,s\,-\,2\,i\,t} \,
 \f{\G(s+\f{k-1}{2}+it)}{\G(1-s+\f{k-1}{2}-it)} \, $$
(\cite[(6.12)]{voronoi}), which has no poles  even in  the larger
region $\Re(s+\f{k-1}{2})>0$.

\section{Choice of functions in the summation
formula}\label{choicesec}

In this section we will describe the test functions $f$ that will
be inserted into (\ref{v}) in order to obtain our eventual
results.  Our goal now is to collect some estimates on $F(x)$ for
the analysis of the righthand side of (\ref{v}) in
\propref{rhsboundprop}. At this stage it is probably helpful to
list which of our variables are considered fixed, and which we
will make estimates in terms of. The parameters $q$, $\l_1,\l_2,$
and $\l_3$ are all considered fixed. At times we will need to
introduce some finite parameters indexed by $\sigma$, $K$, $M$, or
$N$, for example to shift contour integrals or integrate by parts;
this amount will always be bounded in terms of the fixed
parameters $q$ and $\{\l_j\}$. The dependence on these latter
parameters -- ultimately traceable back -- will not be explicitly
mentioned, though it is possible and obviously cumbersome to do
so. The estimates on $F(x)$ in this section and the next mainly
involve the quantities $x$ and $Y$ (a non-negative parameter); the
most important aspect of the bounds on $F(x)$  is their dependence
in $Y$ for $Y \ge 1$.

In order to use (\ref{MFandMf}) it is necessary that $f$ (and
hence $F$) be of parity $\eta\in \Z/2\Z$; we shall accordingly
describe  choices of $f$ for both parities. To make the notation
uniform and convenient, we will from now on regard the parameter
$\d_3$ as an element of $\{0,1\}$, not just of $\Z/2\Z$. Let
$\omega \in \Z/2\Z$ be an arbitrary parity parameter, and fix a
smooth function $\phi_0$ of parity $\eta+\omega\in\Z/2\Z$ with
support in the interval $(-1,1)$.
 From $\phi_0$ we will
define a number of auxiliary functions in terms of the
non-negative parameter $Y$. First will be
\begin{equation}\label{phi1def}
    \phi(x) \ \ = \ \ \left\{%
\begin{array}{ll}
    \ \ \ \ \ \  \phi_0(x-Y) \ + \ (-1)^{\,\omega} \, \phi_0(x+Y)\,,
     & \d_3\equiv 0 \pmod 2 \\
     \\
    \f{1}{2\pi i} \, \left[ \, \phi_0'(x-Y) \ + \ (-1)^{\,\omega}
     \, \phi_0'(x+Y)
        \,  \right], & \d_3\equiv 1 \pmod 2 \, . \\
\end{array}%
\right.
\end{equation}
We let $f(x)\, =\, |x|^{\l_3\,} \sg(x)^{\d_3}\,
\widehat{\phi}(x)$, so that
\begin{equation}\label{fdef}
   f(x) \ \ = \ \ \left\{%
\begin{array}{ll}
   \ \ \ \  2 \, |x|^{\l_3+\d_3} \,\cos(2\pi Y x) \,
    \widehat{\phi}_0(x)\, , & \omega\equiv 0 \pmod 2 \\
    \\
      -2 \, i \, |x|^{\l_3+\d_3}  \ \sin(2\pi Y x) \,
      \widehat{\phi}_0(x) \, , & \omega \equiv 1 \pmod 2 \, . \\
\end{array}%
\right.
\end{equation}
Clearly $\phi$ and $\widehat{\phi}$ are Schwartz functions, so $f$
is admissible in (\ref{v}) and has parity $\eta$.  We have now
\begin{equation}\label{ringoffire}
    M_\eta f (s)   \ \ = \ \  M_{\d_3+\eta}\,\widehat{\phi}(s+\l_3)
    \ \ = \ \ (-1)^{\d_3+\eta} \, G_{\d_3+\eta}(s+\l_3) \,
    M_{\d_3+\eta\,}\phi(1-s-\l_3)
\end{equation}
by (\ref{mellhatg}), and
\begin{equation}\label{boynamedsue}
M_\eta F(s-1) \ \  = \ \  G_{\d_1+\eta}(s-\l_1)
                        \ G_{\d_2+\eta}(s-\l_2)
                        \ M_{\d_3+\eta\,}\phi(s-\l_3) \, ,
\end{equation}
by (\ref{MFandMf}) and (\ref{grecip}). This last expression is
holomorphic in $\Re{s}\ge \f 12$ by (\ref{holominsg1}), because
the signed Mellin transform $M_{\d_3+\eta\,}\phi(s)$ of a Schwartz
function
 can only have poles where
$G_{\d_3+\eta}(s)$ does (see \cite[(3.31)]{inforder}).  Moreover,
$M_\eta F(s)$ decays rapidly in vertical strips.  We may therefore
calculate $F(x)$ using the Mellin inversion formula, e.g.
\begin{equation}\label{MFcountertheline} \gathered
    \f{4 \, \pi \,  i  \, F(x)}{|x| \, \sg(x)^\eta} \ \ = \ \
    \int_{\Re{s} \, =  \, \sigma  \, \ge  \, \f 12}  M_{\eta\,} F(s-1)  \ |x|^{-s}  \ ds \ \ =
 \qquad \qquad \qquad \qquad \\   = \ \ \int_{\Re{s} \, = \,
 \sigma \, \ge \, \f 12}G_{\d_1+\eta}(s-\l_1)
                        \ G_{\d_2+\eta}(s-\l_2)
                        \ M_{\d_3+\eta\,}\phi(s-\l_3) \ |x|^{-s}\ ds \,
                        ,
\endgathered
\end{equation}
or
\begin{equation}\label{MFcountertheline1}\gathered
    \f{4 \, \pi  \, i  \, F(x \,Y)}{|x\, Y| \sg(x)^\eta}\,Y^{\l_3} \ \
    = \qquad\qquad\qquad\qquad\qquad\qquad\qquad\qquad
    \qquad\qquad\qquad\\
     = \ \ \int_{\Re{s} \, =  \, \sigma  \, \ge  \, \f 12}G_{\d_1+\eta}(s-\l_1)
                        \ G_{\d_2+\eta}(s-\l_2)
                        \ M_{\d_3+\eta\,}\phi_1(s-\l_3) \ |x|^{-s}\ ds \, .
\endgathered
\end{equation}
In this last expression, which will be useful for $Y$ large, we
have introduced the function $\phi_1(x) = \phi(Yx)$, which has
$M_{\d_3+\eta\,}\phi(s) = Y^s M_{\d_3+\eta\,}\phi_1(s)$.

\begin{lem}\label{Fxsmall} For $Y\le 1$ we have the uniform
estimate
\begin{equation}\label{Fxsmallpowerbd}
    F(x)  \ \ \ll  \ \ |x|^{\,-N}  \hbox{~~~~~for any real number~~}
N \, \ge \,  \hbox{$-\f 12$}\ ,
    \end{equation}
where the implied constant depends continuously on $N$.
\end{lem}
{\bf Proof:}  We have just remarked above that the integrand in
(\ref{MFcountertheline}) is holomorphic for $\Re{s}=\sigma\ge \f
12$.
 For $\sigma$ in this range
\begin{equation}\label{deliasgone}
    F(x) \ \  \ll \ \  |x|^{1-\sigma} \
    \int_{\R} \left| G_{\d_1+\eta}(\sigma+it-\l_1)
    G_{\d_2+\eta}(\sigma+it-\l_2)  M_{\d_3+\eta}\phi(\sigma+it-\l_3)
      \right| \, dt \, .
\end{equation}
The lemma will follow with $\sigma=N+1$ once we show the integral
in (\ref{deliasgone}) is bounded independently of $Y\le 1$.  To
estimate the function $G$ along vertical lines, we  use the
asymptotic
\begin{equation}\label{stir}
  \left|  G(\sigma+it) \right|  \ \ \sim  \ \ \(\f {|t|}{2\pi}
    \)^{\sigma-1/2}\ , \  \ \ \ \ t \, \rightarrow \, \infty\ ,
\end{equation}
which is a direct consequence of Stirling's formula applied to
definition (\ref{Gdcossin}). Bounds on
$M_{\d_3+\eta\,}\phi(\sigma+it)$ can be obtained from
$$M_{\d_3+\eta\,}\phi(\sigma+it) \ \ = \ \
\int_{\R} \phi(x)\,|x|^{\sigma+it-1}\,\sg(x)^{\d_3+\eta}\,dx \ \le
\ \int_{\R} |\phi(x)|\, |x|^{\sigma-1}\,dx\, .
$$
Because $\phi$ is supported in $(-2,2)$ and is bounded
independently of $Y\le 1$, this last integral is uniformly bounded
in $Y\le 1$ with a continuous dependence on $\sigma\ge \f 12$. The
same holds true when $\phi$ is replaced by any of its derivatives
$\phi^{(K)}$, so $M_{\d_3+\eta+K}\phi^{(K)}(\sigma+it+K)$ is
uniformly bounded in $t$ and $Y$ for $Y\le 1$. Integration by
parts $K$ times then shows
$$M_{\d_3+\eta\,} \phi(\sigma+it)  \ \ \ll \ \  \, |t|^{-K} \
|\, M_{\d_3+\eta+K\,} \phi^{(K)}(\sigma+it+K) \, |  \ \ \ll \ \
|t|^{-K}\, ,$$ again uniformly for $Y\le 1$. Consequently the
integral in (\ref{deliasgone}) converges rapidly and is bounded
independently of $Y\le 1$, with a continuous dependence on
$\sigma$. \bx

The situation for $Y\ge 1$ is more complicated.
 A helpful difference is that $\phi$ vanishes in a neighborhood
 of the origin when $Y\ge  1$, making its Mellin transform
 entire.  We first state a lemma about the Mellin transform's
 dependence on $Y$:

\begin{lem}\label{mphi1}
For any real numbers $Y \ge 1$, $N\ge 0$, and $\sigma$,
\begin{equation}\label{Fxbig3}
    M_{\d_3+\eta\,}\phi_1(\sigma+it-\l_3) \ \ \ll \ \ |t|^{-N}
    |\,M_{\d_3+\eta+N\,}\phi_1^{(N)}(\sigma+it-\l_3+N)\,| \ \  \ll  \ \ \f{1}{Y}\(\f
    Y{|t|}\)^{N}\, .
\end{equation}
\end{lem}
\noindent Here the implied constants are independent of $Y$ and
$t$, and depend continuously on $\sigma$ and $N$.

 {\bf Proof:}
The first inequality comes directly from integration by parts.  We
shall prove the second
 when $N\ge 0$ is an integer; it extends to reals by
 interpolation.  We have
\begin{equation}
\label{Fxbig2} \gathered
    M_{\d_3+\eta+N\,} \phi_1^{(N)}(\sigma+ it) \ \ = \ \
    2 \int_0^\infty \,\phi_1^{(N)}(x) \ x^{\,\sigma\,+\,it\,-\,1} \ dx
    \qquad\qquad\qquad\qquad\qquad
    \\
   \qquad\qquad\qquad \ll \ \ Y^N  \int_0^\infty  \,
   \left|\, \phi_0^{(N+\d_3)}(Y(x-1)) \, \right|
     \, |x|^{\sigma-1} \, dx \, .
\endgathered
\end{equation}
The derivatives of $\phi_0$ are bounded by an absolute constant,
and furthermore the integrand above is supported in the interval
$(1- c \, Y\i,1+ c\,Y\i)$ for some absolute constant $0<c<1$. So
(\ref{Fxbig2}) is bounded by
$O\(Y^{N}\int_{1-c\,Y\i}^{1+c\,Y\i}x^{\sigma-1} dx\)$  $=
O(Y^{N-1})$ for $Y\ge 1$, where the implied constant depends
continuously on $\sigma$. This establishes the second inequality
in (\ref{Fxbig3}). \bx

The arguments used earlier to bound $F(x)$ when $Y\le 1$
generalize to the case $Y\ge 1$ as well, but with an inadequate
$Y$-dependence.  To improve upon \lemref{Fxsmall} for $|x|\le Y$,
we shift the contour further to the left and estimate the
contribution from the poles, rather than merely  a contour
integral positioned just to their right. (The estimates for $|x|
\ge Y\ge 1$ will be given in \propref{xmed} at the end of the next
section.)

\begin{lem}\label{Fxbigsmall}
    For $Y\ge 1$ and $x\le Y$ we have the bound
    \begin{equation}\label{Fxbigsmall1}
    F(x) \ \ \ll \ \ |x|^{1/2}\  Y^{\,-\,{1}/{2}\,-\,\Re \l_3} \, .
\end{equation}
\end{lem}

{\bf Proof:} Suppose momentarily that the poles of
$G_{\d_1+\eta}(s-\l_1)$ and $G_{\d_2+\eta}(s-\l_2)$ do not
overlap.  Shifting the contour in (\ref{MFcountertheline1}) to
$\sigma$ sufficiently negative and avoiding the poles, we obtain
the expression
\begin{equation}\label{Fxbigsmall2}
    \f{4 \, \pi \, i \, F(x \, Y)}{|x|\sg(x)^\eta} \,Y^{\l_3-1} \ \ = \ \
    \sum_{j \, = \, 1,2} \, \sum_{0 \,\le \, k \, < \,  \Re{\l_j}-\sigma}  c_{k,j}\,
    |x|^{\,k-\l_j} \,
    M_{\d_3+\eta\,}\phi_1(\l_j-\l_3-k) \ + \  R\, ,
\end{equation}
where
\begin{equation}\label{Fxbigsmall3}\gathered
    R \ \ \ll \qquad\qquad\qquad\qquad\qquad\qquad\qquad\qquad
    \qquad\qquad\qquad\qquad\qquad\qquad\qquad\qquad \\
     |x|^{-\sigma}
    \int_{\Re{s}\,=\,\sigma\,\ll\, 0} \left| G_{\d_1+\eta}(\sigma+it-\l_1)
    \,
    G_{\d_2+\eta}(\sigma+it-\l_2) \,  M_{\d_3+\eta\,}\phi_1(\sigma+it-\l_3)
      \right| \, dt \\ \ll \ \
\f{|x|^{-\sigma}}{Y} \int_{\R} \(|t|+1\)^{2\sigma+\Re{\l_3}-1}
\,dt \ \ \ll \ \ \f{|x|^{-\sigma}}{Y}\ .
\endgathered
\end{equation}  Here
we have used the fact $M_{\d_3+\eta}\phi_1(s)$ is entire,
(\ref{normassumpt}), (\ref{stir}), and (\ref{Fxbig3}) with $N=0$;
also the implied constants in (\ref{Fxbigsmall3}) depend
continuously on $\sigma$. Finally the $c_{k,j}$ are constants
coming from the residues of $G_{\d_1+\eta}(s-\l_1)
G_{\d_2+\eta}(s-\l_2)$ at the points $s=\l_j-k$. Another
application of (\ref{Fxbig3}) bounds the $M_{\d_3+\eta}\phi_1$
factor in the sum on the righthand side of (\ref{Fxbigsmall2})  by
$O(1/Y)$ as well.

If the poles of $G_{\d_1+\eta}(s-\l_1)$ and
$G_{\d_2+\eta}(s-\l_2)$ in fact do overlap, then
(\ref{Fxbigsmall2}) remains correct provided an additional factor
of $(\log |x|)^{j-1}$ is included. Bounding the righthand side of
(\ref{Fxbigsmall2}) therefore gives the estimate
\begin{equation}\label{Fxbigsmall4}
    F(x\,Y) \, Y^{\l_3}  \ \ \ll_\e  \ \
    |x|^{\,1-\Re\l_1-\e}\ + \ |x|^{\,1-\Re\l_2-\e}\ \ , \ \  \ \
    |x| \ \le \  1
\end{equation}
for any $\e>0$.  These exponents are both greater than $\f 12$ for
$\e$ small, thanks to (\ref{orderoflambda}), and so the righthand
side of (\ref{Fxbigsmall4}) is bounded by $O(|x|^{1/2})$, for $x
\le 1$,
 proving (\ref{Fxbigsmall1}).
 \bx

\section{A substitute for stationary phase}\label{subssec}

To bound $F(x)$ from the integral (\ref{vvv}) one could attempt to
use stationary phase. We instead find it more convenient to apply
a device of \cite[p.~33-]{nara} to the Mellin transform $MF$ of
$F$ instead. This allows us to express the transformed $F$ in
terms of an asymptotic series of Fourier transforms, which in
practice are often simpler to estimate. Let us use the notation
$f(s) \approxeq g(s)$ if, for any integer $M\ge 0$, $f(s)/g(s)$
has an asymptotic expansion of the form $1+ c_1 s\i+\cdots + c_M
s^{-M}+O(|s|^{-M-1})$ for large values of $s$
 in any vertical strip of finite width.
In this notation Stirling's formula reads $\G(s)\, \approxeq \,
\sqrt{2\pi}\,e^{-s} \,s^{\,s-1/2\,}$.  Consequently,
\begin{equation}\label{stirfry}
    \f{\G(s+a)}{\G(s+b)} \ \  \approxeq  \ \ s^{\,a-b} \,.
\end{equation}

Though we only need a special case, we will state the following
lemma  in enough generality that it can be applied to arbitrary
$L$-functions.  Indeed, the ratio of $\G$-factors in the
functional equation of any $L$-function can always be written in
the form of the lefthand side of (\ref{johnnyc}) below (see
\cite[\S 6]{voronoi}, for example).  The method here can often be
applied instead of stationary phase on $\R^n$ to give asymptotic
expansions of the transformed functions in general Voronoi-style
summation formulas (e.g.~\cite{expos}) and approximate functional
equations, in terms of the ordinary, one-variable Fourier
transform.

\begin{lem}\label{jcash}
For any $\(\mu_1,\ldots,\mu_n\)\in \C^n$ with mean $\bar\mu$, and
$(\e_1,\ldots,\e_n)\in (\Z/2\Z)^n$, we have
\begin{equation}\label{johnnyc}
    \prod_{j\,=\,1}^n \, G_{\e_j}(s-\mu_j) \ \  \ \approxeq  \ \ \
   \sum_{\g\,=\,0}^1  \,  C_\g \ n^{-n \, s}  \, \  G_{\g}\(n \, s  -
   n \,\bar\mu
   + \f{1-n}{2}\)
 \end{equation}
for explicitly computable constants $C_0$ and $C_1$.
\end{lem}
{\bf Proof:} To simplify the notation, add $\bar\mu$ to $s$ and
denote the sum $\sum_{j=1}^n \e_j$ as simply $\e$.  We can thereby
assume that $\bar{\mu}=0$. The lefthand side of (\ref{johnnyc})
can be expressed using (\ref{Gdcossin}) as
\begin{equation}\label{johnnyc2}\gathered
    \prod_{j=1}^n \ (2\pi)^{\,-s-\mu_j}\ \G(s-\mu_j) \ \left[
      e\(\f{s-\mu_j}{4}\)  \ + \  (-1)^{\e_j} \ e\(\f{\mu_j-s}{4}\)
\right]  \ \ \ \ \  \qquad
\\
\qquad\qquad\qquad\approxeq  \ \ (2\pi)^{-\,n\,s} \left[\, e\(\f{n
s}{4}\)\, + \, (-1)^\e \, e\(-\f{ns}{4}\)\, \right]
\,\prod_{j=1}^n  \, \G(s-\mu_j) \,.
\endgathered
\end{equation}
We now use (\ref{stirfry}) and the identity
\begin{equation}\label{triplication}
    \prod_{j=0}^{n-1} \, \G\(s+ \f{j}{n}\) \ \ = \ \ (2\pi)^{\f{n-1}{2}}
    \, \
     n^{\,\f12 \,-\,n\,s}  \
    \G(ns)
\end{equation}
to rewrite (\ref{johnnyc2}) as
\begin{equation}\label{johnnyc3}
    \aligned
\approxeq \ \ &   (2\pi)^{-\,n\,s} \left[\, e\(\f{n s}{4}\)\, + \,
(-1)^\e \, e\(-\f{ns}{4}\)\, \right]  \ \prod_{j=1}^n  \,
\G\(s+\f{j-1}{n}\) s^{-\,\mu_j\,-\,(j-1)/n}
\\
\approxeq  \ \ & (2\pi)^{-\,n\,s+ (n-1)/2} \left[\, e\(\f{n
s}{4}\)\, + \, (-1)^\e \, e\(-\f{ns}{4}\)\, \right]  \,
s^{\,(1-n)/2} \  n^{\,1/2\,-\, n\, s} \  \G(ns)
\\ \approxeq  \ \ &
(2\pi)^{\,-n\, s+(n-1)/2}  \left[\, e\(\f{n s}{4}\)\, + \, (-1)^\e
\, e\(-\f{ns}{4}\)\, \right] \,n^{\,n/2\,-\,n\,s} \ \G\(n\, s
\,+\,\f{1-n}{2}\)  .
   \endaligned
\end{equation}
It is  clear from (\ref{Gdcossin}) that $G_0(s)\pm G_1(s)$ span
all linear combinations of $(2\pi)^{-s} \,\G(s)\,e(\pm \f s4)$, so
this last expression is indeed a linear combination of
$n^{-\,n\,s}\,G_\g(n\,s+\f{1-n}{2})$, $\g=0,1$, as (\ref{johnnyc})
asserts.
 \bx

We next remark that the same argument used in the last step of
(\ref{johnnyc3}) shows that
\begin{equation}\label{gears}
   \f{ G_\d(s+1)}{s}  \ \ = \ \ (2\pi)^{\,-s-1}\ \G(s) \
    \left[ \,
e\(\f{s+1}{4}\) \   +  \ (-1)^\d \  e\(-\f{s+1}{4}\)
    \, \right]
\end{equation}
is also a linear combination of $G_0(s)$ and $G_1(s)$.  That means
that the higher terms in the asymptotic expansion in
(\ref{johnnyc}) can also be written in terms of linear
combinations of the $G_\d$'s, with shifted arguments.  We shall
now apply this specifically to the product of $G_\d$'s in
(\ref{boynamedsue}):
\begin{equation}\label{suegs}
    G_{\d_1+\eta}(s-\l_1)\, G_{\d_2+\eta}(s-\l_2) \ \ \approxeq \ \
\sum_{\g\,=\,0}^1 \,C_\g \, 2^{-2s}\,G_\g(2 \,s
\,-\,\l_1\,-\,\l_2\,-\,\hbox{$\f 12$})\, ,
\end{equation}
 and return to bounding $F(x)$ in the regime $|x| \ge Y \ge 1$.
We can use (\ref{normassumpt}) and (\ref{suegs})  to restate
(\ref{boynamedsue}) as
\begin{equation}\label{MFevenapprox2}\gathered
    M_\eta F(s-1) \ \
\qquad \qquad \qquad \qquad \qquad \qquad \qquad \qquad \qquad
\qquad \qquad \qquad \qquad
     \\    = \ \ \scriptstyle{   \sum_{\g=0}^1 \(
 \sum_{j=0}^{M-1} C_{\g,j} \, 2^{-2s} \,
    G_\g(2s+\l_3-\f 12 -j) \, + \,
    O(2^{-2s} G_\g(2s+\l_3-\f 12 -M)) \) \,  M_{\d_3+\eta}\phi(s-\l_3)}.
    \endgathered
\end{equation}
The error term represented by the $O$-notation here comes from the
asymptotic expansion; the implied constants of course depend only
on $\l_1$ and $\l_2$, which we consider fixed. We will take $M$ to
be a large positive integer, and evaluate $F(x)$ using
(\ref{MFevenapprox2}) in the contour integral representation
(\ref{MFcountertheline1}) along $\Re{s}=\sigma$, where
$2\sigma+\Re{\l_3}-M =\e$, an arbitrarily small positive real
number.  Recall the remark at the beginning of \secref{choicesec}
that $M$ will be bounded in terms of $q$ and the $\{\l_j\}$.
Changing the value of the constants $C_{\g,j}$, we may write
\begin{equation}\label{bigchain1} \gathered
    Y^{\l_3}\sg(x)^\eta \, \f{F(x\,Y)}{|x|\,Y} \ \ =
\qquad \qquad \qquad \qquad \qquad \qquad \qquad\qquad \qquad
\qquad \qquad  \\ = \ \
    \sum_{\g\,=\,0}^1
\(    \sum_{j\,=\,0}^{M-1} C_{\g,j} \int_{\Re{s}\,=\,\sigma}
|4x|^{-s}\,
    G_\g(2s+\l_3-\hbox{$\f 12$}-j)\, M_{\d_3+\eta\,}\phi_1(s-\l_3) \, ds \right. \\
\left.    +  \ O\(|4x|^{-\sigma}  \int_{\Re{s}\,=\,\sigma}
     \left|G_\g (2\sigma+2it+\l_3-\hbox{$\f 12$} -M) \,
       M_{\d_3+\eta\,}\phi_1(\sigma+it-\l_3) \right|
  \,    dt  \) \) \,
\endgathered
\end{equation}

The sum of these error terms, which we denote $R$, can be bounded
by (\ref{stir}) and \lemref{mphi1} as
\begin{equation}\label{bigchain2}\gathered
    R \ \ \ll
    \qquad\qquad\qquad\qquad\qquad\qquad\qquad\qquad
    \qquad\qquad\qquad\qquad\qquad\qquad\qquad\\
     |x|^{-\sigma} \( \int_0^1 \f{dt}{Y} \,+\,\int_{1}^Y \f{
    t^{\,2\,\sigma\,+\,\Re{\l_3} \,-\, M \,-\, 1 }}{Y} \, dt  \, + \,
    \int_{Y}^\infty Y^{N-1}\,
    t^{\,2\,\sigma\,+\,\Re{\l_3} \,-\, M \,-\, 1 \,-\, N} \, dt \) \\
   = \ \ O\( |x|^{-\sigma} \, Y^{\,2\,\sigma\,+\,\Re{\l_3}\,-\,M\,-\,1} \) \ \
   = \ \ O\(\,\left|\f{x}{Y^2}\right|^{-\,\sigma}
   Y^{\,\Re{\l_3} \,-\, M \,-\,1}\,\)
\endgathered
\end{equation}
for any large $N\ge 0$ (recall $2\sigma+\Re\l_3-M=\e>0$). Changing
variables in (\ref{bigchain1}) gives an  expression -- again with
different constants $C_{\g,j}$ -- of the form
\begin{equation}\label{bigchain3}
    Y^{\l_3-1}\sg(x)^\eta\,\f{F(x\,Y)}{|x|} \, -  \, R \ \ = \ \
    \sum_{\g\,=\,0}^1
    \sum_{j\,=\,0}^{M-1}\ C_{\g,j}  \
    \left| x \right|^{\,(\l_3-j-1/2)/2} \
    \widehat{\psi}_{\,\g,j}\(2|x|^{1/2}\)
  \,,
\end{equation}
where $$M_{\g\,}\psi_{\,\g,j\,}(s) \ \   = \ \
M_{\d_3+\eta\,}\phi_1\(\f{1-s-3\l_3+\f 12 +j}{2}\),$$  or in other
words,
$$\psi_{\,\g,j\,}(x) \ \ = \  \ 2 \, \sg(x)^{\,\d_3\,+\,\eta\,+\,\g}
\ \, |x|^{\,3\,\l_3\,-\,j\,-\,\f 32} \ \,
 \phi_1\(1/x^2\)\, .$$  Recall from (\ref{phi1def}) that each
$\psi_{\g,j}$ is smooth, and supported in a neighborhoods of width
$O(Y\i)$ about $\pm 1$.   In addition, a straightforward
calculation writing $\psi_{\g,j}$ in the form
$\sg(x)^{\d'\,}|x|^{A\,}\phi(\f{Y}{x^2})$ shows that the $K$-th
derivative of $\psi_{\g,j}$ is bounded by $O(Y^K)$. We conclude
that
\begin{equation}\label{psihatparts}
\widehat{\psi}(r)  \ \ \ll \ \  r^{-K} \ ||\,\psi^{(K)\,}||_1
 \ \ \ll \ \ r^{-K} \, Y^{K-1}\,
\end{equation}
 for each of our functions $\psi=\psi_{\g,j}$, and  any integer
$K\ge 0$.

Inserting (\ref{bigchain2}) and (\ref{psihatparts}) into
(\ref{bigchain3}), we obtain the following bound for $F(x)$ when
$|x|\ge 1$:
\begin{equation}\label{beeneverywhere}
    \f{Y^{\l_3-1}}{|x|}\,F(x\,Y) \ \ \ll \ \
 Y^{K\,-\,1}\
     |x|^{-\,\f{K}{2}\,+\,\f{\Re{\l_3}}{2}\,-\,\f 14}
  \ \, + \, \    \left|\f {x}{Y^2}
    \right|^{-\sigma}\,Y^{\,\Re{\l_3}-M-1}  \ ,
\end{equation}
\begin{equation}\label{tuscaloosa}
    F(x\,Y) \ \ \ll \ \  Y^{K\,-\,\Re{\l_3}}\ |x|^{\f 34\,-\,\f
    K2\,+\,
    \f{\Re{\l_3}}{2}}
    \ \, + \, \   Y^{\,2\,\sigma\,-\,M} \ |x|^{1-\sigma}
   \ ,
\end{equation}
or
\begin{equation}\label{tuscaloosb}
    F(x) \ \ \ll \ \  Y^{\,\f {3K}{2}\,-\, \f 34 \, -\, \f{3\,\Re{\l_3}}{2}}
    \ |x|^{\f{3}{4}\,-\,\f K2\,+\,\f{\Re{\l_3}}{2}} \
    \, + \, \ Y^{\,3\,\sigma\,-\,M\,-\,1}\  |x|^{1-\sigma}
   \ .
\end{equation}
Recalling our choice that $2\sigma+\Re{\l_3}-M =\e$, where $\e>0$
is arbitrarily small, we may deduce our final estimates on $F(x)$:
 \begin{prop}\label{xmed}
Let $Y\ge 1$.

(a)  If $\ Y \, \le \,  |x| \,  \le \,  Y^3$ then for any
$\e\,>\,0$
\begin{equation}\label{xmed2}
    F(x) \ \ \ll_{\,\e} \ \ Y^{\,\e\,-\,\Re{\l_3}} \ \left| \f xY
    \right|^{\,\f 34 \,+\,\f{\Re{\l_3}}{2}}  .
\end{equation}

 (b)  If   $\ |x| \,  \ge
\, Y^3$ then for any $N\,>\,0$
\begin{equation}\label{xmed1}
    F(x) \ \ \ll_{\,N}  \ \ Y^{\,3/2}\ \,
    \left| \f{x}{Y^3}\right|^{-\,N}\,.
\end{equation}
\end{prop}
{\bf Proof:}  These both follow directly from (\ref{tuscaloosb}).
For part (a), set $K=0$ to handle the first term, and note that
$\sigma$ is large in the second.  To settle part (b) it suffices
to prove the prove (\ref{xmed1}) for $N$ large, which is
straightforward because $\sigma$, $M$, and $K$ may be taken to be
large.\bx

{\bf Remark.} The method used here can be used to obtain more
precise information about the asymptotic behavior of $F(x)$.  In
particular, the fact that $\psi_{\g,j}(x)$ is concentrated near
$x=\pm 1$ allows one to understand the oscillatory behavior of
$F(x)$ as well (see \cite{nara2}, where this is explored in much
more detail for summation formulas connected to Dedekind zeta
functions).

\section{Proof of \thmref{mainthm}}\label{pfsec}

In this section we prove \thmref{mainthm} by inserting our choice
of $f$ into the summation formula (\ref{v}).  First we need to
specify some of the parameters in that formula. Let $T\ge 1$,
$a=-p_k$,  and $c=q_k$, where $p_j/q_j$ are the continued fraction
approximants to $\a$, and $k$ is chosen such that
\begin{equation}\label{Tandcontfrac}
    q_k^2 \ \ \le \ \  T \ \ \ \le \ \ q_{k+1}^2.
\end{equation}
We set $Y=T\,|\a+\f ac|\ge 0 $ so that $\a= \pm \f YT-\f ac$, and
hence
\begin{equation}\label{cfrac}
  \f {1}{q_k \, q_{k+1}} \ \ \ll \ \
     \f YT  \  \ = \ \  \left| \,  \a-\f{p_k}{q_k}  \,
     \right| \ \ \ll \ \ \f {1}{q_k \,
     q_{k+1}}
\end{equation}
by the standard properties of continued fractions (see, for
example, \cite[p.~47]{baker}).

The following proposition applies our bounds on $F$ to the
righthand side of (\ref{v}) in \thmref{vsf}; afterwards we will
conclude \thmref{mainthm} by a standard analysis of the lefthand
side.

\begin{prop}\label{rhsboundprop}
With the choice of $f$ given in (\ref{phi1def}-\ref{fdef}), the
righthand side of (\ref{v}) is $O_\e(T^{3/4+\e})$, independent of
$\a$.
\end{prop}
{\bf Proof:}
 First, note that the GCD of the
 parameters  of the
 Kloosterman sum in (\ref{v}) is bounded by $GCD(q\bar{a},qc/d)\le q$,
  which we consider
 to be fixed.  This Kloosterman sum is therefore bounded by
$O_\e((qc/d)^{1/2+\e})=O_\e((c/d)^{1/2+\e})$, for any $\e>0$,
according to Weil's bound (the implied constant here of course
depends on $\e$). That means the righthand side of (\ref{v}) is
bounded by
\begin{equation}\label{rhsv1}
    \sum_{d|c} \left| \f cd\right|^{3/2+\e\,} \sum_{n \neq 0}
     \f{|a_{n,d}|}{|n|} \left|F\(\f{nd^2}{c^3q} T\)
     \right|.
\end{equation}
Now we will use our bounds on $F$ to bound this expression. First,
let us settle the case of $Y\le 1$, which is simpler. Here we
break up (\ref{rhsv1}) as the sum of
\begin{equation}\label{rhsv2}
    \sum_{d|c} \left| \f cd\right|^{3/2+\e}
    \sum_{n \, \neq \,  0\,, \ |n| \, \le \, X}
     \f{|a_{n,d}|}{|n|} \, \left|F\(\f{n}{X} \)
     \right|
\end{equation}
and
\begin{equation}\label{rhsv3}
    \sum_{d|c} \left| \f cd\right|^{3/2+\e}
    \sum_{ |n| \, >  \, X}
     \f{|a_{n,d}|}{|n|}  \, \left|F\(\f{n}{X} \)
     \right|\, ,
\end{equation}
 where for convenience we have set $X=\f{c^3 q}{d^2T}$. We will
use the  bound of $F(x) \ll x^{-N}$ from  \lemref{Fxsmall} here:
for (\ref{rhsv2}) we take $N=0$, and for (\ref{rhsv3}) we take
$N$ arbitrarily large.
 Then
the sums over $n$ are therefore bounded by
\begin{equation}\label{rhsv4a}
\sum_{0  \, \neq  \, |n| \,  \le \,  X} \, |a_{n,d}|   \, |n|^{-1}
\ \ \  \ \hbox{and} \ \ \ \ \sum_{|n| \,
> \,  X} \, |a_{n,d}| \, |X|^{N} \,|n|^{-1-N}\,,
\end{equation}
 respectively, with the remaining value of $N$ a large positive
integer.  Let us assume momentarily that $X\ge 1$.  The partial
summation identity
\begin{equation}\label{partsum}
    \sum_{n\,=\,1}^K\,a_n \,b_n \ \ = \ \
    \sum_{n\,=\,1}^{K\,-\,1}\, A_n\,(b_n\,-\,b_{n+1}) \,+ \,
    A_K\,b_K \  \ ,  \ \ \ A_n \ = \ \sum_{k\,=\,1}^n \,a_n
\end{equation}
can now be used to replace the coefficients $|a_{n,d}|$ in
(\ref{rhsv4a}) by their average size of $O(d)$ from
(\ref{ramave}). One then bounds the two sums by $O(d \log
(X+1))=O_\e(d X^\e)$ and
\begin{equation}\label{intapr}
\ll \ \ d \, \int_X^\infty \f{dn}{n} \( \f nX      \)^{-N} \ \ = \
\ d \, \int_1^\infty n^{-N-1}\,dn \ \ = \ \  O(d) \, ,
\end{equation}
 respectively.
  Thus the righthand side of
(\ref{v}) is
\begin{equation}\label{endoffirstpart}
 \ll_{\,\e}  \ \ |c|^{\,3/2\,+\,\e\,}
\,\sum_{d|c}d^{\,-\,1/2\,-\,\e} \, X^\e \ \ \ll_{\,\e} \ \
|c|^{\,3/2\,+\,\e} \, X^\e  \ \ = \  \ O_\e(T^{\,3/4\,+\,\e})
\end{equation}
 by (\ref{Tandcontfrac}), and the fact that
$\#\{d|c\}=O_\e(c^{\,\e})$.  This has been done subject to the
assumption that $X\ge 1$, but actually the argument simplifies if
$X<1$ because the first sum in (\ref{rhsv4a}) has no terms and
(\ref{intapr}) is taken over a shorter range.

Now we turn to the case where $Y\ge 1$, which is more involved. We
now break the sum over $n$ in (\ref{rhsv1}) into three ranges:
\begin{eqnarray}
    \sum_{n \,  \neq \, 0}
     \f{|a_{n,d}|}{|n|} \left|F\(\f{n}{X}\)
     \right|  & \le & \label{rhsv4b} \\
     & &   \sum_{0 \, \neq \, |n| \,  < \,  X\,Y }
     \f{|a_{n,d}|}{|n|} \left|F\(\f{n}{X} \)
     \right|  \label{low} \\
     & +&   \sum_{ X \, Y \ \le \  |n| \  \le  \ X\,Y^3 }
     \f{|a_{n,d}|}{|n|} \left|F\(\f{n}{X} \)
     \right|  \label{med} \\
     & +&   \sum_{|n| \ > \ X\,Y^3 }
     \f{|a_{n,d}|}{|n|} \left|F\(\f{n}{X} \)
     \right| \, , \label{high}
\end{eqnarray}
with again $X=\f{c^3 q}{d^2T}$.  We will again make the assumption
that $XY\ge 1$; otherwise the analysis is  simpler as it was just
above in the argument for $Y\le 1$.  For (\ref{low}) we use the
estimate $F(n/X) \ll Y^{-\Re\l_3}$ from \lemref{Fxbigsmall}. After
again using partial summation to replace $|a_{n,d}|$ by $O(d)$,
this results in the bound
\begin{equation}\label{lowbd}
\gathered
 \ll_{\,\e} \ \   d\ Y^{-\Re\l_3}
\sum_{|n| \, \le \, XY}
    \f 1{|n|}  \ \ \ll \ \ d\  Y^{-\Re{\l_3}} \ \log(XY+1) \ \ \ll_\e
    d\ (XY)^\e \,
.
\endgathered
\end{equation}
for (\ref{low}) (recall $\Re{\l_3}\ge 0$).

For the remaining pieces (\ref{med}-\ref{high}), we turn to
\propref{xmed}. The bound (\ref{xmed2}) allows us to bound
(\ref{med}) by
\begin{equation}\label{medbd} \ll_{\,\e}
  \ \ X^{-D} \
 Y^{\,\e\,-\,\Re{\l_3} \,-\,D}
  \   \sum_{ | n |\,  \le \, XY^3}
 |a_{n,d}|\, |n|^{\,D\,-\,1} \, ,
\end{equation}
where $D=\f{3}{4}+\f{\Re\l_3}{2} >0$.  We then estimate
(\ref{medbd}) again by partial summation, and find it is
\begin{equation}\label{medbd1}
  \ll_\e \ \ d \  X^{-D} \ Y^{\,\e-\Re{\l_3}-D} \ (XY^3)^{D} \ \
    = \ \ d \ Y^{\,\e\,-\,\Re{\l_3}\,+\,2D} \ \ = \ \ d \ Y^{\,3/2\,+\,\e}\,.
\end{equation}

Finally, for (\ref{high}) we use the bound of (\ref{xmed1}),
namely $F(n/X) \ll Y^{\,3/2} |\f{n}{XY^3}|^{-N}$ for any $N$
large. Again after removing the $|a_{n,d}|$ by partial summation,
(\ref{high}) is bounded by
\begin{equation}\label{highbd}
    \ll \ \ d \ Y^{3/2} \, |XY^3|^{N} \,
    \sum_{n\ge XY^3}  \, |n|^{-N-1}  \ \ \ll \ \  d \ Y^{3/2} \, .
\end{equation}

All told  (\ref{rhsv4b}) is $O(d\, Y^{3/2+\e})$, plus the
negligible term (\ref{lowbd}), which is $O_\e(d\, T^\e\, Y^\e)$.
 The final contribution of the sum over $d$ in (\ref{rhsv1})
is again $O_\e(c^{\e})$ just as immediately after
(\ref{endoffirstpart}), so the righthand side of (\ref{v}) is
bounded by
$$O_\e(T^{\,\e}\,(c\,Y)^{3/2+\e}) \ \ = \ \
O_\e\(T^\e \(\f{T}{q_{k+1}}\)^{3/2+\e\,}\) \ \ = \ \
O_{\e\,}(T^{3/4+\e})\, ,$$ because of (\ref{Tandcontfrac}) and
(\ref{cfrac}). \bx

 We have just bounded the righthand side of (\ref{v}).  By
 taking linear combinations
of the lefthand side for the functions (\ref{fdef}) for
$\eta,\omega \in \Z/2\Z$, one obtains the result
\begin{equation}\label{punch1}
    \sum_{n\neq 0} a_{n,q}  \, |n|^{\l_3+\d_3}  \, e(n\a) \,
    \widehat{\phi}_0\(\f nT \) \ \  =  \ \ O_\e \( T^{\l_3+\d_3+3/4+\e}\) \, ,
\end{equation}
uniformly in $\a\in \R$, for any  smooth function $\phi_0$ with
support in $(-1,1)$.  Here in (\ref{punch1}) the implied constant
depends also on $\phi_0$.  By rescaling $\phi_0$ with the
parameter $T$, (\ref{punch1}) remains valid for any smooth
function $\phi_0$ of compact support.

\begin{lem}\label{dkern}
For any  $g \in C^1(\R)$, the $L^1$ norm of the kernel
\begin{equation}\label{perkern}
    D_{g,N}(x) \ \ = \ \ \sum_{n \, = \,  1}^N \ e(nx) \,
g\(\f{n}{N}\)
\end{equation}
taken over $\R/\Z$ satisfies $$||D_{g,N}||_1  \ \ = \ \  O(\log
N)\, ,$$ where the implied constant depends on $g$.
\end{lem}
{\bf Proof:} By the  partial summation formula (\ref{partsum}),
\begin{equation}\label{dkern2}
    D_{g,N}(x) \ \ = \ \ \sum_{n\,=\,1}^{N-1} D_n(x) \, \(g(\f{n}{N})
    -g(\f{n+1}{N})\) \ + \ D_N(x)\,g(\f{N}{N})\,,
\end{equation}
where $D_n(x) = \sum_{k\le n}e(kx)=
\f{e((n+1)x)\,-\,e(x)}{e(x)\,-\,1}$ is essentially the classical
Dirichlet kernel, and has $||D_n||_1 =O(\log n)$. The $L^1$ norm
of $D_{g,N}$ is thus bounded by
\begin{equation}\label{dkern3}
     ||D_{g,N}||_1 \ \ \le \ \ \sum_{n=0}^{N-1}\, ||D_n||_1 \,
     \left|g(\f{n}{N})
    -g(\f{n+1}{N})\right|  \ +  \ ||D_N||_1 \,
    \left|g(\f{N}{N})\right|\,.
\end{equation}
As  $g\(\f{n}{N}\) - g\(\f{n+1}{N}\)=O(1/N)$, (\ref{dkern3}) is
bounded by $O(N\f{\log N}{N}) + O(\log N) = O(\log N)$.\bx

{\bf Proof of \thmref{mainthm}:} We have shown (\ref{punch1})
holds for any  smooth function $\phi_0$ with compact support.
Choose $\phi_0$ so that $\widehat{\phi}_0$ is non-zero on $[0,1]$
(that such a function exists can be seen simply by rescaling).
Letting $T=N$ and $g=1/\widehat{\phi}_0(x)$, convolve the lefthand
side of (\ref{punch1}) over $\R/\Z$ against the kernel  $D_{g,N}$
from (\ref{perkern}).  The uniform upper bound in (\ref{punch1})
and $L^1$-norm estimate from \lemref{dkern} provides us the
following estimate:
\begin{equation}\label{punch2}
    \sum_{n  \,  \le   \,
    T} \  a_{n,q}  \  n^{\l_3 \, + \, \d_3} \   e(n \, \a) \ \
   = \ \  O_{\e}\( T^{ \, \l_3 \, + \, \d_3 \, + \, 3/4 \, + \, \e} \) \, .
\end{equation}
The power of $n$ may be then removed using the partial summation
formula  (\ref{partsum}), proving \thmref{mainthm}.\bx

\section{Proof of \thmref{implic}}
\label{moments}

In this section we prove \thmref{implic}, which gives a bound on
the second moment of an $L$-function based on the amount of
cancellation present in additive twists of its coefficients.  Our
assumption in \thmref{implic} is that the $L$-function $L(s) =
\sum_{n=1}^\infty a_n n^{-s}$, assumed to not be the Riemann
$\zeta$-function, is the standard $L$-function of a cusp form on
$GL(m)$ over $\Q$.  We now quickly review their basic analytic
properties (see for example \cite{godjac} and
\cites{gelmil,rudsar,gelsha,jacquetindia,bumprs}).
   The coefficients $a_n$ are of polynomial growth
 and $L(s)$ satisfies a
functional equation of the form
\begin{equation}\label{fe}
    L(s) \ \ = \ \ \omega \, N^{1/2-s} \, {\mathcal G}(1-s) \,
    \tilde{L}(1-s)\,,
\end{equation}
where $N\ge 1$ is the ``conductor,'' $\omega$ is a complex number
of modulus $1$, and ${\mathcal G}(s)$ is a ratio of $\G$-factors.
The latter is customarily written in a variety of styles, though
this choice is inessential here; for a later stage in this
argument it will be convenient  to write
\begin{equation}\label{curlyG}
    {\mathcal G}(s) \ \ = \ \ \prod_{j\,=\,1}^m\,G_{\d_j}(s+\l_j)\,,
\end{equation}
where $\d_j\in \Z/2\Z$, $\l_j \in \C$, and $G_\d(s)$ is the
function defined in (\ref{Gdcossin}).\footnote{The parameters
$\{(\d_j,\l_j)\}$ can be viewed as principal series embedding
parameters for the representation of $GL(m,\R)$ associated to the
cusp form. In fact they have already made an appearance in
\thmref{vsf}; for a fuller discussion of the case $m=3$, see
\cite{voronoi}.} Finally the dual $L$-function in (\ref{fe}) is
defined by $\tilde{L}(s) = \sum_{n=1}^\infty
\overline{a_n}n^{-s}$, and both $L(s)$ and $\tilde{L}(s)$ are
entire and of finite order, except for the excluded case of $m=1$
and $L(s)=\zeta(s)$. Furthermore, $L(s)$ vanishes at certain
points where the $\G$-factors have poles; all that we will utilize
is that there exists an index $1\le k \le m$ with $\Re{\l_k}\le 0$
such that $G_{\d_k}(s-\l_k)L(s)$ is also entire. For shorthand we
denote the pair $(\d_k,\l_k)$ as $(\d,\nu)$.

At the heart of the connection between cancellation bounds and the
second moment is the classical method of Titchmarsh \cite[p.
165]{Titch}, which uses Parseval's identity for the Mellin
transform (our conventions here are carried over from
\secref{vsec}). Let $\phi(x) \in C_c^\infty(\R)$ have parity
$\d\in\Z/2\Z$. The function
\begin{equation}\label{conv1}
    f(x) \ \ = \ \ \sum_{n\,\neq\, 0} \,
    a_{|n|}\,\sg(n)^{\d}\ \widehat\phi(nx)\ |nx|^{-\nu}
\end{equation}
also has parity $\d$, and its signed Mellin transform is
\begin{equation}\label{conv2}
\aligned
    M_\d f(s) \ \ & = \ \ \sum_{n\neq 0} \,
    a_{|n|}\,\sg(n)^\d\,\int_{\R}\widehat{\phi}(nx)\,|nx|^{-\nu}
    \,|x|^{s-1}\,\sg(x)^\d\,dx \\
& = \ \ 2\,L(s)\,M_{\d\,}\widehat\phi(s-\nu) \\
& = \ \ 2 \,(-1)^\d\, L(s) \, G_\d(s-\nu) \,
M_{\d\,}\phi(1-s+\nu)\, \qquad\qquad\text{(cf. (\ref{mellhatg})).}
\endaligned
\end{equation}
Parseval's identity
\begin{equation}\label{parsf}
 \int_\R |f(x)|^2\,dx \ \ = \ \ \f{1}{4\pi} \, \int_\R\,
    |M_\d f({\scriptstyle \f 12} + it)|^2\,dt\,
\end{equation}
relates the second moment of $L(\f 12 + it)$ to the $L^2$ norm of
$f$ as follows.  Because of  (\ref{stir}), one has
\begin{equation}\label{conv3}
\aligned
    |M_\d f({\scriptstyle \f 12} + it)| \ \ & = \ \ 2 \, |L({\scriptstyle \f 12}
     + it )| \,|G_\d({\scriptstyle \f 12} +it-\nu)|\,
    |M_{\d}\phi({\scriptstyle \f 12} - it+\nu)|\,
\\
& \sim \ \ 2 \,({\scriptstyle \f {|t|}{2\pi}})^{-\Re\nu}
|L({\scriptstyle \f 12}
     + it )| \,
    |M_{\d}\phi({\scriptstyle \f 12} - it+\nu)|\,.
\endaligned
\end{equation}
We shall now pick $\phi$ more specifically in order to bound the
second moment through (\ref{parsf}).   The main idea is to ensure
that $M_{\d\,} \phi(\f 12 - it + \nu)$, $t\in\R$, approximates the
characteristic function of $[-T,T]$, for  $T$ is a large real
parameter. Let $\phi_0\ge 0 $ be an even smooth function supported
in the interval $(-\f 12,\f 12)$, and let
\begin{equation}\label{convphi1} \phi(x) \ \ = \ \ \phi_0(T(x-1))
\ + \ (-1)^\d\,\phi_0(T(x+1))\,,
\end{equation} so that
\begin{equation}\label{convphi2}
    \widehat\phi(r) \ \ = \ \ \left[ e(-r) \, + \,
     (-1)^\d \,e(r)        \right]\, T\i \,
    \widehat{\phi_0}\(\f rT \)\,
\end{equation}
and
\begin{equation}\label{convphi3}
    f(x) \ \ = \ \ \sum_{n\neq 0} \, a_{|n|}
    \,\sg(n)^\d\,|nx|^{-\nu}\,
    [e(-nx)\,+\,(-1)^\d\,e(nx)]\ T\i \,
    \widehat{\phi_0}\(\f{nx}{T}\)\,.
\end{equation}
For $T$ large, the function $\phi_0(T(x-1))$ is concentrated near
$x=1$ and has mass on the order of $T\i$.  It is straightforward
to choose $\phi_0$ such that $M_{\d\,}\phi(\f 12 - it + \nu)$ is
nonzero in the range $t\in [-cT,cT]$ for some $c>0$. By rescaling
$\phi_0$ if necessary, one may further ensure $|M_{\d\,}\phi(\f 12
- it + \nu)|\gg T\i$ for $t\in[-T,T]$; as a result of
(\ref{parsf}),
\begin{multline}\label{convphi4}
   T^{-2}\ \int_{-T}^{T} \, |t|^{\,-\,2 \,\Re \nu}\ |L({\scriptstyle \f
    12}+it)|^2\ dt \ \ \ll_{\,\e} \\ \ \ \int_{|x| \, \le \,  T^{1-m-\e}}
    |f(x)|^2\,dx \ + \ \int_{|x| \, \ge \,  T^{1+\e}}
    |f(x)|^2\,dx \ + \ \int_{ T^{1-m-\e} \,  < \,  |x| \,  < \, T^{1+\e}}
    |f(x)|^2\,dx\,.
\end{multline}
Since $\phi$ is supported away from the origin for $T$ large, its
Mellin transform is entire.  Therefore the last expression for
$M_\d f(s)$ in (\ref{conv2}) is also entire because of our
assumption on $\nu$.  A standard contour shift and  application of
Stirling's formula (or alternatively the asymptotic analysis
developed for $F$ in \secref{subssec}), produces rapid decay of
$f$ near $0$ and $\infty$ -- enough to make the first two terms on
the righthand side of (\ref{convphi4}) decay rapidly in $T$
because of the ranges so chosen. The expression for $f(x)$ in
(\ref{convphi3}) is a linear combination of  smoothed variants of
$|x|^{-\nu} T\i S(\f T x,x)$, except of course  for the added
presence of the $|n|^{-\nu}$ term. These differences can be
removed by partial summation as at the end of  the last section,
and so our assumption that $S(T,x)=O_\e(T^{\b+\e})$ gives the
bound $f(x) = O_\e(|x|^{-\Re\nu\;}T\i |\f T x|^{\b-\Re\nu+\e})$.
We conclude that the righthand side of (\ref{convphi4}) is
$O_\e(T^{-2\,\Re\nu-1-m+2\,\b \, m+\e})$.
  Recalling that $\Re\nu \le 0$, this implies that
  $\int_{-T}^T|L(\f 12 + it )|^2dt=O_\e(T^{\,2\,\b \, m-m+1+\e})$.

\bx

\section{Period bounds: an analytic analog of the Ramanujan
conjecture}

In this section we establish an equivalence of
 the
folklore cancellation conjecture (\ref{sdef}-\ref{best}) for the
$L$-function coefficients  of a cusp form $\phi$ on
$\G_1(N)\backslash GL(m,\R)$, where
 $\G_1(N)$ denotes the subgroup of $GL(m,\Z)$
consisting of matrices whose last row equals $[0\,0\,\cdots\,0\,1]
\pmod N$. The methods and results in this section are not truly
particular to the subgroups $\G_1(N)$ themselves, but this family
is a very canonical one to study because it captures every cusp
form on $GL(m)$: namely, every adelic automorphic representation
 has a
vector which is (left-)invariant under $\G_1(N)$
\cite{jpss-conductor}. So there is essentially no loss of
generality entailed by this restriction.

Let us first introduce the period alluded to above, which first
originated in the construction of the standard $L$-function on
$GL(m)$ by Jacquet, Piatetski-Shapiro, and Shalika
\cite{maryland,jpss}
 (see also
\cites{Cogdell,bumprs,jacquetindia} as references for this
section). Let $P$ denote the standard $(2,1,1,\ldots,1)$ parabolic
of $GL(m)$, so that its unipotent subgroup $N$ consists of all
upper triangular unit matrices which have zero as their second
entry in the first row (blank entries are zero):
\begin{equation}\label{pandn}
    P=\left\{\(\begin{smallmatrix}
      \star & \star & \star & \star & \star & \star \\
      \star & \star & \star & \star & \star & \star \\
         &   & \star & \star & \star & \star \\
       &    &   &  \star & \star & \star \\
       &    &     &   & {}_{\ddots} & \star \\
       &  &  &  &  & \star \\
    \end{smallmatrix}        \)  \ \in \ GL(m) \right\}  \ , \ \
     N=\left\{\(\begin{smallmatrix}
      1& 0 & \star & \star & \star & \star \\
       & 1 & \star & \star & \star & \star \\
       &  & 1 & \star & \star & \star \\
       &    &   &  1 & \star & \star \\
       &    &     &   & {}_{\ddots} & \star \\
       &  &  &  &  & 1\\
    \end{smallmatrix}        \)  \ \in \ GL(m) \right\} \,.
\end{equation}
  Let $\psi$
denote the standard additive character of unit upper triangular
matrices, which maps a matrix $n$ to $e^{\,2\,\pi\, i \,s(n)}$,
$s(n)$ being the sum of the entries of $n$ lying one position
above the diagonal. Clearly $\psi$ is invariant under $N(\Z)$.
 The
period under consideration is,
\begin{equation}\label{psperiod}
    V(g) \ \ = \ \ \int_{N(\Z)\backslash N(\R)} \phi(ng) \,
    \overline{\psi(n)}\,dn\,.
\end{equation}
Of course our notation suppresses the implicit dependence of $V$
on $\phi$.

We can now state the main result of this section.  For $m=2$ this
result is well-known, though the history (especially of the
implication (b) $\Longrightarrow$ (a)) is somewhat muddled -- the
first complete proof  we are aware of is in \cite{hafner}.

\begin{thm}\label{perthm}
Let $0\neq\phi\in C^\infty(\G_1(N)\backslash GL(m,\R))$ be a cusp
form and $\f 12 \le \b<1$. Then the following two statements are
equivalent:

(a)  The $L$-function coefficients $a_{n}$, $n\ge 1$, satisfy the
cancellation bound
\begin{equation}\label{thmabd}
    \sum_{n\,=\,1}^T \, a_n  \, e(n x) \ \ = \ \
   O_\e( T^{\,\b\,+\,\e})
\end{equation}
uniformly in $x$ for any $\e>0$;

(b)  The period $V$ satisfies the bound
\begin{equation}\label{thmbbd}
   V\(\begin{smallmatrix}
  y & x & & \\ {} & 1 & & \\ & & \ddots & \\ & & & 1 \end{smallmatrix} \) \ \
   = \ \ O_\e\(\,y^{\,\f{m-1}{2}\,-\,\b\,-\,\e}\,\) \ , \ \ \ y \, > \,
   0\ ,
\end{equation}
 uniformly in $x$ for any $\e>0$.
\end{thm}

{\bf Remarks:}  1)  The reason we have termed this an ``analytic
analog of the Ramanujan conjecture'' is that the conjectured
optimal bound in (\ref{thmbbd}) with $\b=\f 12$  is reminiscent of
the following classical situation.  Let $\phi(z) \, = \,
\sum_{n\ge 1}\,c_n\, e(nz)$ be the Fourier expansion of a
classical holomorphic form of weight $k$. The trivial bound on the
$n$th-coefficient
\begin{equation}\label{cnrmaper}
    c_n  \ \ =  \ \ e^{ \, 2 \, \pi \, n \, y}  \,
    \int_{x\,=\,0}^1 \,\phi(x+iy) \,e(-n x) \, dx \ \  \le
    \ \  e^{\,2 \, \pi \,
n \,  y} \,  \int_{x\, = \, 0}^1 \, |\phi(x+iy)| \, dx\ \
\end{equation}
is obtained by invoking the bound $\phi(x+iy)=O(y^{-k/2})$ and
taking $y$ to be of order $1/n$: $c_n = O(n^{k/2})$.  This
estimate is on the order of $\sqrt{n}$ short of the truth  of
$c_n=O(n^{(k-1)/2})$ predicted by the Ramanujan conjecture (in
this case a theorem of Deligne \cite{Deligne}). The reason for
this loss of $\sqrt{n}$
 in (\ref{cnrmaper})  is
that we have used absolute values and forfeited any cancellation
from the oscillation of this period integral. A similar phenomenon
likely happens in (\ref{thmbbd}), for bounding (\ref{psperiod})
trivially via absolute values presumably gives an estimate which
is off by some power of $n$. The analogy with the Ramanujan
conjecture is only meant in this analytic sense and is not meant
to have any algebraic connotation.

2) Note that the period in part (b), like $\phi$, decays rapidly
as $y\rightarrow \infty$.  So the (\ref{thmbbd}) is only an issue
for $y$ small.

3) Since different vectors in the same representation space share
common $L$-function coefficients, assertion (b) is either true for
all or none of the nonzero smooth vectors in the irreducible
subrepresentation of $L^2(\G_1(N)\backslash GL(m,\R))$ generated
by right translates of $\phi$.

4) Finally, the reason we have focused on the range $\f 12 \le \b
< 1$ is because in practice this is only interesting situation
($\b=1$ being trivial, and $\b=\f 12$ conjectured to be optimal).

One of the advantages of taking $\phi$ to be invariant under
$\G_1(N)$ is the Fourier expansion (see \cite{shalika,psexpn} and
\cite[(2.1.6)]{bumprs})
\begin{multline}\label{psexpn}
    \phi(g) \ \ = \ \ \\
    \sum_{n_1,\,\ldots,\,n_{m-1}\,\ge\,1}\,\sum_{\g \, \in \,
    \G_\infty\backslash
    GL(m-1,\Z)}\,\f{a_{\,n_1,\,\ldots\,,\,n_{m-1}}}{\prod_{j\,=\,1}^{m-1}
    n_j^{\ j\,(m-j)/2}}\, W\( \(\begin{smallmatrix}  n_1 & & &
& \\  & n_2 & & &
\\  & & \ddots & &
\\
& & & n_{m-1} &
\\ & & & & 1  \end{smallmatrix}\)\ttwo{\g}{}{}1  g \) \, ,
\end{multline}
where $\G_\infty$ refers to the subgroup of unit upper triangular
matrices in $GL(m-1,\Z)$.  Here $W(g)$ is the (archimedean)
Whittaker function, formed from the same type of integral as
(\ref{psperiod}), but with $N$ replaced by the maximal unipotent
subgroup $N_0=\{$all unit upper triangular matrices$\}$. The
coefficients of the standard $L$-function in this
notation\footnote{The identification of the coefficients in
(\ref{psformula}) with those of the standard $L$-function (as
opposed to the contragredient $L$-function) is somewhat arbitrary
and  not completely universal; in fact we used a different
convention above in (\ref{L-functions}).  The difference is of no
essential consequence here, and was  introduced purely as a matter
of convenience. } are $a_n=a_{n,1,1,\ldots,1}$.
 Automorphic
representations always satisfy the transformation law
$\phi(gh)=(-1)^\d\phi(g)$ for some $\d\in \Z/2\Z$, where
$h=\(\begin{smallmatrix} -1 & & &  & \\  & 1 & & &
\\  & & 1 & &
\\
& & & \ddots &
\\ & & & & 1  \end{smallmatrix}\)$; the Whittaker
functions naturally inherit this right-transformation property as
well.
 Formulas
(\ref{psperiod}) and (\ref{psexpn}) then give the following
 expression for $V(g)$ in terms of the $L$-function
coefficients:
\begin{equation}\label{psformula}
    V(g) \ \ = \ \ \sum_{n\,\neq\, 0} \ \f{a_{|n|}}{|n|^{(m-1)/2}}\
W\(\(\begin{smallmatrix} n & & &  & \\  & 1 & & &  \\  & & 1 & &
\\
& & & \ddots &
\\ & & & & 1  \end{smallmatrix}  \)g\)\,.
\end{equation}
  The integral representation of the standard $L$-function on
$GL(m)$ by Jacquet, Piatetski-Shapiro, and Shalika (see
\cite{jpss,jpss2,jsblue,Cogdell,jacquetindia,bumprs}) uses the
signed Mellin transform
\begin{multline}\label{lfnint}
    I(s,\phi) \ \ = \ \ \int_\R\,V\(\begin{smallmatrix}  y & & &  & \\  & 1 & & &
\\  & & 1 & &
\\
& & & \ddots &
\\ & & & & 1  \end{smallmatrix}\)\,|y|^{\,s-\f{m-1}{2}-1}\,\sg(y)^\d\,dy  \ \ =
\\ = \ \
2\,J(s,W)\,\sum_{n\,=\,1}^\infty\,a_n \,n^{-s} \ \ = \ \
2\,J(s,W)\,L(s,\phi)\, ,
\end{multline}
where
\begin{equation}\label{Jint}
   J(s,W) \ \ = \ \  \int_\R  W\(\begin{smallmatrix}  y & & &  & \\  & 1 & & &
\\  & & 1 & &
\\
& & & \ddots &
\\ & & & & 1  \end{smallmatrix}\)\,|y|^{\,s-\f{m-1}{2}-1}\,\sg(y)^\d\,dy\,.
\end{equation}
It is known by a result of Jacquet and Shalika
\cite{jsblue,jacshalikavol} that $J(s,W)$ can have poles only
where the $\G$-factors in the functional equation of $L(s,\phi)$
do.  Because of their unitary bound \cite{Jac-Sha}, none of these
poles have $\Re{s}\ge \f 12$ (see \cite[(2.5) and
(5.18)]{rudsar}).  Consequently,
\begin{equation}\label{whitstate}
    W\(\begin{smallmatrix}  y & & &  & \\  & 1 & & &
\\  & & 1 & &
\\
& & & \ddots &
\\ & & & & 1  \end{smallmatrix}\) \ \ = \ \ O(y^{\,m/2-1}) \ , \ \ \
y \, > \, 0 \, ,
\end{equation}
because the lefthand side possesses an asymptotic expansion as
$y\rightarrow 0$ (in fact much more is known -- see for example
\cite{flato,wallach}.)  This bound, and the entirety of $J(s,W)$
in $\Re{s}\ge 1/2$, will be used in the proof which follows.

{\bf Proof of \thmref{perthm}:}

 The key link between parts (a) and
(b) is formula (\ref{psformula}).  First we treat the case  (a)
$\Longrightarrow$ (b) using partial summation. The archimedean
Whittaker functions of smooth vectors are themselves smooth, decay
rapidly, and -- thanks to (\ref{whitstate}) --  have only mild
growth near the origin. Let
\begin{equation}\label{oneparam}
w(y)\ \ = \ \ W\(\begin{smallmatrix}  y & & &  & \\  & 1 & & &
\\  & & 1 & &
\\
& & & \ddots &
\\ & & & & 1  \end{smallmatrix}\),
\end{equation} so that the formula for $V\(\begin{smallmatrix}
  y & x & & \\ {} & 1 & & \\ & & \ddots & \\ & & & 1 \end{smallmatrix} \)$
  (denoted more succinctly as just  $V\(\begin{smallmatrix}
  y & x \\ {} & 1 \end{smallmatrix} \)$) in (\ref{psformula})
  reads
\begin{equation}\label{psform2}
    V\(\begin{smallmatrix}
  y & x \\ {} & 1 \end{smallmatrix} \) \ \ = \ \ y^{\,(m-1)/2}\
  \sum_{n\,\neq\,0}
  \, a_n\, \sg(n)^\d \,  e(nx)\ |ny|^{(1-m)/2}\ w(|ny|)\ , \ \  y \, > \, 0\,.
\end{equation}
Letting $f(t)=|t|^{(1-m)/2\,}w(|t|)$ and keeping in mind that $y$
may be assumed to be small (Remark 2), we bound (\ref{psform2}) by
partial summation:
\begin{equation}\label{psform3}
       V\(\begin{smallmatrix}
  y & x \\ {} & 1 \end{smallmatrix} \)  \ \ \ll_\e \ \
  y^{\,(m-1)/2}\  \sum_{n\,=\,1}^\infty \ n^{\,\b\,+\,\e}\,y
  \ |f'(ny)|\  \, ,  \ \ \  y \, > \, 0\,.
\end{equation}
The Whittaker function $w(t)$ (or, more accurately, the
restriction of the Whittaker function to the one parameter
subgroup in (\ref{oneparam})) has rapid decay as $t\rightarrow
\infty$, which means that in order to achieve (\ref{thmbbd}) we
need just to establish the bound
\begin{equation}\label{psform4}
    \sum_{n\,=\,1}^{T} \,n^{\,\b\,+\,\e} \,|f'(ny)| \ \ = \ \ O_\e(y^{-1-\b-\e})\,,
\end{equation}
 where $T$ is on the order of $y^{-1-\e'}$
for $\e'$ very close to 0, say $\e'=\e^2$.  This itself follows
from knowing that
\begin{equation}\label{psform5}
\gathered
 f'(t) \ \ = \ \  {\scriptstyle \f{1-m}{2}}\  t^{-(m+1)/2} \, w(t) \  \, + \,  \
t^{(1-m)/2} \, w'(t) \ \ = \ \ O(t^A) \ \ , \ \ \ t\ > \ 0 \\
\qquad \qquad \qquad \qquad \qquad \qquad \qquad  \qquad
\qquad\text{~~for some~~} A\,\ge\,-\,\b-1\,
 ,
\endgathered
\end{equation}
 for then the sum in (\ref{psform4}) is bounded by
$O_\e(y^A\sum_{n=1}^T n^{\b+A+\e}) =O_\e(y^{-1-\b-\e})$ (for $\e'$
sufficiently small). The condition (\ref{psform5}) in turn follows
from
\begin{equation}\label{psform6}
    w(t) \ , \ \ t\,w'(t)  \ \ = \ \ O(t^{\,m/2-1}) \, ,
\end{equation}
where we have used the assumption that $1/2\le\b<1$.  The last
assertion is literally  (\ref{whitstate}) for $w(t)$ itself, and
in fact also for $t w'(t)$ too, which is the Whittaker function of
the derivative $Ad(e_{11})\phi$ of $\phi$, $e_{11}$ being the
matrix which is zero everywhere except for the entry $1$ in its
first position.  This completes the proof of (a) $\Longrightarrow$
(b).

The reverse implication (b) $\Longrightarrow$ (a) can be proven
using the technique of Hafner \cite{hafner}. Namely, the rapid
decay as $y\rightarrow\infty$ and the bound as $y\rightarrow 0$ of
(\ref{thmbbd}) gives an analytic continuation of the Mellin
transform of
\begin{equation}\label{haf1}
    \int_\R \,V\(\begin{smallmatrix}
  y & x \\ {} & 1 \end{smallmatrix} \) \,
  |y|^{\,s-\f{m-1}{2}-1}\,\sg(y)^\d\,dy \ \ = \ \
  J(s,W)\ \sum_{n\neq 0}^\infty \, a_n \, e(nx) \, |n|^{-s}
\end{equation}
(cf. (\ref{lfnint})) to $\Re{s}
> \b$. In this range $J(s,W)$ is holomorphic by the remark above
(\ref{whitstate}).   For real $s>\beta$ such that $J(s,W)\neq 0$,
(\ref{haf1}) gives the bound
\begin{equation}\label{haf2}
    \sum_{n\,\neq\,0}\,a_{n}\,e(nx)\,  |n|^{-s}  \ \ =  \ \ O(1) \, ,
\end{equation}
uniformly in $x$, and with the implied constant depending on $s$.
In fact the sum on the lefthand side is a continuous function of
$x$ because of dominated convergence applied to the integral in
(\ref{haf1}). This in turn implies (a) by \cite[Prop.
3.7]{regularity}.

\bx

 The implication (b) $\Longrightarrow$ (a) can alternatively  be proven using the techniques in
\cite{flato}, or instead by the method of Voronoi summation in
\cite{voronoi}.

\vspace{2 cm}{\noindent Stephen D. Miller \newline Department of
Mathematics
\newline Hill Center-Busch Campus
\newline Rutgers, The State University of New Jersey
\newline 110 Frelinghuysen Rd
\newline  Piscataway, NJ 08854-8019
\newline {\tt miller@math.rutgers.edu}
}


\begin{thebibliography}{99}


\bib{baker}{book}{
    author={Baker, Alan},
     title={A concise introduction to the theory of numbers},
 publisher={Cambridge University Press},
     place={Cambridge},
      date={1984},
     pages={xiii+95},
      isbn={0-521-24383-1},
      isbn={0-521-28654-9},
    review={MR 86f:11001},
}



\bib{Bump}{book}{
     author={Bump, Daniel},
      title={Automorphic forms on ${\rm GL}(3,{\bf R})$},
     series={Lecture Notes in Mathematics},
     volume={1083},
  publisher={Springer-Verlag},
      place={Berlin},
       date={1984},
      pages={xi+184},
}

\bib{bumprs}{article}{
    author={Bump, Daniel},
     title={The Rankin-Selberg method: a survey},
 booktitle={Number theory, trace formulas and discrete groups (Oslo, 1987)},
     pages={49\ndash 109},
 publisher={Academic Press},
     place={Boston, MA},
      date={1989},
    review={MR 90m:11079},
}

\bib{nara}{article}{
    author={Chandrasekharan, K.},
    author={Narasimhan, Raghavan},
     title={The approximate functional equation for a class of
            zeta-functions},
   journal={Math. Ann.},
    volume={152},
      date={1963},
     pages={30\ndash 64},
    review={MR 27 \#3605},
}

\bib{nara2}{article}{
    author={Chandrasekharan, K.},
    author={Narasimhan, Raghavan},
     title={Exponential sums associated with algebraic number fields},
   journal={Comment. Math. Helv.},
    volume={54},
      date={1979},
    number={4},
     pages={523\ndash 561},
      issn={0010-2571},
    review={MR 81j:10057},
}

\bib{chowla}{article}
{author={Chowla, S.D.}, title={Some problems of diophantine
approximation (I)}, journal={Math. Zeitschrift}, volume={33},
year={1935}, pages={544\ndash 563}}



\bib{Cogdell}{book}{
author ={Cogdell, J. W.}, title={$L$-functions and Converse
Theorems for $GL(n)$}, series = {IAS/Park City Lecture Notes},
place={Park City, Utah}, date={2002}, note = {\newline
\url{http://www.math.okstate.edu/~cogdell/}}}

\bib{Deligne}{article}{
    author={Deligne, Pierre},
     title={La conjecture de Weil. I},
  language={French},
   journal={Inst. Hautes \'Etudes Sci. Publ. Math.},
    number={43},
      date={1974},
     pages={273\ndash 307},
    review={MR 49 \#5013},
}

\bib{sarnak}{article}{
    author={Epstein, Charles},
    author={Hafner, James Lee},
    author={Sarnak, Peter},
     title={Zeros of $L$-functions attached to Maass forms},
   journal={Math. Z.},
    volume={190},
      date={1985},
    number={1},
     pages={113\ndash 128},
      issn={0025-5874},
    review={MR 86m:11034},
}


\bib{erdos}{article}{
    author={Erd{\"o}s, P.},
     title={Some remarks on Diophantine approximations},
   journal={J. Indian Math. Soc. (N.S.)},
    volume={12},
      date={1948},
     pages={67\ndash 74},
    review={MR 10,513b},
}

\bib{gelmil}{article}{
    author={Gelbart, Stephen S.},
    author={Miller, Stephen D.},
     title={Riemann's zeta function and beyond},
   journal={Bull. Amer. Math. Soc. (N.S.)},
    volume={41},
      date={2004},
    number={1},
     pages={59\ndash 112 (electronic)},
      issn={0273-0979},
    review={2 015 450},
}



\bib{gelsha}{book}{
    author={Gelbart, Stephen},
    author={Shahidi, Freydoon},
     title={Analytic properties of automorphic $L$-functions},
    series={Perspectives in Mathematics},
    volume={6},
 publisher={Academic Press Inc.},
     place={Boston, MA},
      date={1988},
     pages={viii+131},
      isbn={0-12-279175-4},
    review={MR 89f:11077},
}


\bib{godjac}{book}{
    author={Godement, Roger},
    author={Jacquet, Herv{\'e}},
     title={Zeta functions of simple algebras},
      note={Lecture Notes in Mathematics, Vol. 260},
 publisher={Springer-Verlag},
     place={Berlin},
      date={1972},
     pages={ix+188},
    review={MR 49 \#7241},
}

\bib{hafner}{article}{
    author={Hafner, James Lee},
     title={Some remarks on odd Maass wave forms (and a correction to:
            ``Zeros of $L$-functions attached to Maass forms'' [Math.\ Z.
            {\bf 190} (1985), no.\ 1, 113--128; MR 86m:11034] by Hafner, C.
            Epstein and P. Sarnak)},
   journal={Math. Z.},
    volume={196},
      date={1987},
    number={1},
     pages={129\ndash 132},
      issn={0025-5874},
    review={MR 89a:11056},
}


\bib{HL}{article}{author={G. Hardy},author={J. Littlewood},title=
{Some problems of Diophantine approximation},journal={Acta
Math.}, volume={37}, year={1914},pages={193\ndash 238}}

\bib{iwaniec}{book}{
    author={Iwaniec, Henryk},
     title={Topics in classical automorphic forms},
    series={Graduate Studies in Mathematics},
    volume={17},
 publisher={American Mathematical Society},
     place={Providence, RI},
      date={1997},
     pages={xii+259},
      isbn={0-8218-0777-3},
    review={MR 98e:11051},
}


\bib{jacquetindia}{article}{
    author={Jacquet, Herv{\'e}},
     title={Dirichlet series for the group ${\rm GL}(n)$},
 booktitle={Automorphic forms, representation theory and arithmetic (Bombay,
            1979)},
    series={Tata Inst. Fund. Res. Studies in Math.},
    volume={10},
     pages={155\ndash 163},
 publisher={Tata Inst. Fundamental Res., Bombay},
      date={1981},
    review={MR 83j:10032},
}

\bib{jacshalikavol}{article}{author={Jacquet, H.},title={Integral
Representation of Whittaker functions}, note={to appear in the
Shalika festschrift}}


\bib{maryland}{book}{
    author={Jacquet, H.},
    author={Pyatetskii-Shapiro, I. I.},
    author={Shalika, J. A.},
     title={Construction of cusp forms on ${\rm GL}(3)$},
      note={Department of Mathematics, University of Maryland, Lecture Note,
            No. 16},
 publisher={Department of Mathematics, University of Maryland, College Park,
            Md.},
      date={1975},
     pages={i+85},
    review={MR 58 \#16522},
}


\bib{jpss}{article}{
    author={Jacquet, Herv{\'e}},
    author={Piatetski-Shapiro, Ilja Iosifovitch},
    author={Shalika, Joseph},
     title={Automorphic forms on ${\rm GL}(3)$.},
   journal={Ann. of Math. (2)},
    volume={109},
      date={1979},
    number={1~and~ 2},
     pages={169\ndash 258},
      issn={0003-486X},
    review={MR 80i:10034b},
}

\bib{jpss-conductor}{article}{
    author={Jacquet, H.},
    author={Piatetski-Shapiro, I. I.},
    author={Shalika, J.},
     title={Conducteur des repr\'esentations du groupe lin\'eaire},
  language={French},
   journal={Math. Ann.},
    volume={256},
      date={1981},
    number={2},
     pages={199\ndash 214},
      issn={0025-5831},
    review={MR 83c:22025},
}

\bib{jpss2}{article}{
    author={Jacquet, H.},
    author={Piatetskii-Shapiro, I. I.},
    author={Shalika, J. A.},
     title={Rankin-Selberg convolutions},
   journal={Amer. J. Math.},
    volume={105},
      date={1983},
    number={2},
     pages={367\ndash 464},
      issn={0002-9327},
    review={MR 85g:11044},
}




\bib{Jac-Sha}{article}{
    author={Jacquet, H.},
    author={Shalika, J. A.},
     title={On Euler products and the classification of automorphic
            representations. I},
   journal={Amer. J. Math.},
    volume={103},
      date={1981},
    number={3},
     pages={499\ndash 558},
      issn={0002-9327},
    review={MR 82m:10050a},
}

\bib{jsblue}{article}{
    author={Jacquet, Herv{\'e}},
    author={Shalika, Joseph},
     title={Rankin-Selberg convolutions: Archimedean theory},
 booktitle={Festschrift in honor of I. I. Piatetski-Shapiro on the occasion
            of his sixtieth birthday, Part I (Ramat Aviv, 1989)},
    series={Israel Math. Conf. Proc.},
    volume={2},
     pages={125\ndash 207},
 publisher={Weizmann},
     place={Jerusalem},
      date={1990},
    review={MR 93d:22022},
}



\bib{Kim}{article}{
    author={Kim, Henry H.},
     title={Functoriality for the exterior square of ${\rm GL}\sb 4$ and the
            symmetric fourth of ${\rm GL}\sb 2$},
   journal={J. Amer. Math. Soc.},
    volume={16},
      date={2003},
    number={1},
     pages={139\ndash 183 (electronic)},
      issn={0894-0347},
    note = {With appendix 1 by Dinakar Ramakrishnan and appendix 2 by Kim
            and Peter Sarnak.
            \newline \url{http://www.ams.org/jourcgi/jour-getitem?pii=S0894034702004101}}
}

\bib{expos}{article}{author = {Miller, Stephen D.}, author =
{Schmid, Wilfried}, title={Summation Formulas, from Poisson and
Voronoi to the Present}, booktitle = {Noncommutative Analysis, in
Honor of Jacques Carmona}, series = { Progress in Mathematics},
volume = { 220}, publisher = {Birkh\"auser}, year={ 2003},
note={\newline
\url{http://www.math.rutgers.edu/~sdmiller/voronoi}}}


\bib{inforder}{article}{author = {Miller, Stephen D.}, author =
{Schmid, Wilfried}, title={Distributions and Analytic Continuation
of Dirichlet Series} , journal={Journal of Functional Analysis},
volume={24}, year= {2004}, pages={ 155 \ndash 220}, note =
{\url{http://arxiv.org/abs/math.NT/0403030}}}

\bib{voronoi}{article}{author = {Miller, Stephen D.}, author =
{Schmid, Wilfried}, title={Automorphic Distributions,
$L$-functions, and Voronoi Summation for $GL(3)$} , journal =
{Annals of Mathematics}, note = {to appear.\newline
\url{http://www.math.rutgers.edu/~sdmiller/voronoi}}}

\bib{regularity}{article}{author = {Miller, Stephen D.}, author
= {Schmid, Wilfried}, title={The Highly Oscillatory Behavior of
Automorphic Distributions for $SL(2)$}, journal = {Letters in
Mathematical Physics}, note={to appear.
 \newline \url{http://arxiv.org/abs/math.NT/0402382}}}

\bib{murty}{article}{
    author={Murty, M. Ram},
    author={Sankaranarayanan, A.},
     title={Averages of exponential twists of the Liouville function},
   journal={Forum Math.},
    volume={14},
      date={2002},
    number={2},
     pages={273\ndash 291},
      issn={0933-7741},
    review={MR 2002m:11078},
    note={\url{http://www.degruyter.de/journals/forum/2002/pdf/14_273.pdf}}
}

\bib{psexpn}{article}{
    author={Piatetski-Shapiro, I. I.},
     title={Multiplicity one theorems},
 booktitle={Automorphic forms, representations and $L$-functions (Proc.
            Sympos. Pure Math., Oregon State Univ., Corvallis, Ore., 1977),
            Part 1},
    series={Proc. Sympos. Pure Math., XXXIII},
     pages={209\ndash 212},
 publisher={Amer. Math. Soc.},
     place={Providence, R.I.},
      date={1979},
    review={MR 81m:22027},
}

\bib{rudsar}{article}{
    author={Rudnick, Ze{\'e}v},
    author={Sarnak, Peter},
     title={Zeros of principal $L$-functions and random matrix theory},
      note={A celebration of John F. Nash, Jr.},
   journal={Duke Math. J.},
    volume={81},
      date={1996},
    number={2},
     pages={269\ndash 322},
}

\bib{flato}{article}{
    author={Schmid, Wilfried},
     title={Automorphic distributions for ${\rm SL}(2,\mathbb R)$},
 booktitle={Conf\'erence Mosh\'e Flato 1999, Vol. I (Dijon)},
    series={Math. Phys. Stud.},
    volume={21},
     pages={345\ndash 387},
 publisher={Kluwer Acad. Publ.},
     place={Dordrecht},
      date={2000},
}

\bib{selclass}{article}{
    author={Selberg, Atle},
     title={Old and new conjectures and results about a class of Dirichlet
            series},
 booktitle={Proceedings of the Amalfi Conference on Analytic Number Theory
            (Maiori, 1989)},
     pages={367\ndash 385},
 publisher={Univ. Salerno},
     place={Salerno},
      date={1992},
    review={MR 94f:11085},
}




\bib{shalika}{article}{
    author={Shalika, J. A.},
     title={On the multiplicity of the spectrum of the space of cusp forms
            of ${\rm GL}\sb{n}$},
   journal={Bull. Amer. Math. Soc.},
    volume={79},
      date={1973},
     pages={454\ndash 461},
    review={MR 52 \#5886},
}


\bib{Titch}{book}{
    author={Titchmarsh, E. C.},
     title={The theory of the Riemann zeta-function},
   edition={2},
      note={Edited and with a preface by D. R. Heath-Brown},
 publisher={The Clarendon Press Oxford University Press},
     place={New York},
      date={1986},
     pages={x+412},
      isbn={0-19-853369-1},
    review={MR 88c:11049},
}






\bib{walfisz}{article}
{author={A. Walfisz}, journal={Math. Zeitschrift},
volume={35},title={\"Uber einige trigonometrische Summen. Zweite
Abhandlung}, year={1932}, pages={774\ndash 778}}


\bib{wallach}{article}{
    author={Wallach, Nolan R.},
     title={Asymptotic expansions of generalized matrix entries of
            representations of real reductive groups},
 booktitle={Lie group representations, I (College Park, Md., 1982/1983)},
    series={Lecture Notes in Math.},
    volume={1024},
     pages={287\ndash 369},
 publisher={Springer},
     place={Berlin},
      date={1983},
    review={MR 85g:22029},
}



\end{thebibliography}
\end{document}